\documentclass[reqno, 10pt]{amsart}

\usepackage{graphicx}
\usepackage{amsfonts}

\newtheorem{theorem}{Theorem}[section]                                          
                          
\newtheorem{lemma}[theorem]{Lemma}

\newtheorem{remark}{Remark}[section]

\newcommand{\no}{\noindent}
\newcommand{\nn}{\nonumber}
\newcommand{\oo}{\infty}

\newcommand{\ra}{\rightarrow}

\newcommand{\<}{\langle}

\makeatletter
\@addtoreset{equation}{section}
\makeatother

\renewcommand\thetable{\thesection.\@arabic\c@table}

\bibdata{prob}
\bibstyle{alpha} 

\begin{document}

\title{Evolution of the interfaces in a two dimensional Potts model}

\thanks{Received: 15 August 2006. Revised: 31 March 2007.}

\maketitle

\begin{center}
{GLAUCO VALLE}\footnote{Work supported by CNPq, FAPERJ.}

\medskip

{\footnotesize{UFRJ - Departamento de m\'etodos estat\'{\i}sticos do Instituto de Matem\'atica. \\
Caixa Postal 68530, 21945-970, Rio de Janeiro, Brasil \\
e-mail:  \rm \texttt{glauco.valle@dme.ufrj.br}}}
\end{center}

\begin{abstract}
We investigate the evolution of the random interfaces in a two dimensional Potts model at zero temperature under Glauber dynamics for some particular initial conditions. We prove that under space-time diffusive scaling the shape of the interfaces converges in probability to the solution of a non-linear parabolic equation. This Law of Large Numbers is obtained from the Hydrodynamic limit of a coupling between an exclusion process and an inhomogeneous one dimensional zero range process with asymmetry at the origin. 
\end{abstract}

\section{Introduction} 
\label{sec:potts}

The Potts model is a stochastic process that describes the random evolution of $q$--states spins on a lattice. It furnishes a wide class of interesting models and problems in statistical mechanics. A standard reference on the subject from the physical point of view is the review paper of Fred Wu \cite{wu}. A fundamental point in the study of the Potts model is to understand the evolution of the random interfaces between regions of constant state. By identifying the Potts model to the description of a microscopic phenomenon and imposing on the system a space and time change of scales, we can ask about the convergence of interfaces to some representative macroscopic curve. This is a source that leads to interesting problems in hydrodynamics for physically relevant models of interacting particle systems. 

In this paper we consider a Potts models under very particular initial conditions. However it has an interesting hydrodynamical behavior whose proof has demanding technical difficulties and it is of interest by itself. Moreover the model to be studied has connections to the hydrodynamical behavior of other relevant interacting particle systems as the zero-range and the exclusion processes. These identifications have already been pointed out by Landim, Olla, Volchan in \cite{landimollavolchan} and our study just completes what they have initiated in that paper as we should make clear later. 

To a better understanding of the model, the main results, the connections to other interacting particle systems and the technical difficulties we have to deal with to obtain the hydrodynamical behavior, we go directly to the formal definition of the Potts model.    

We consider a 3-states Potts model at zero temperature under Glauber dynamics. It is described by a spin system on $\Omega_{sp}=\{-1,0,1\}^{\mathbb{Z}^2}$, whose generator acting on cylinder functions is given by
$$
(\mathcal{L}_{sp} f)(\sigma) = \frac{1}{2} \sum_{x\in \mathbb{Z}^2} \sum_{j=-1}^{1} \mathbf{1} \{\mathcal{H}(\sigma^{x,j}) - \mathcal{H}(\sigma) \le 0\} [f(\sigma^{x,j}) - f(\sigma)],
$$
for all $\sigma \in \Omega_{sp}$, where $\sigma^{x,j}$ denotes the configuration
$$
(\sigma^{x,j})(y) = \left\{ \begin{array}{ll}
\sigma(y) & \textrm{for } y\neq x \\
j & \textrm{for } y=x \, ,
\end{array} \right.
$$
and $\mathcal{H}$ is the Hamiltonian defined formally on $\Omega_{sp}$ by
$$
\mathcal{H} (\sigma) = \sum_{x \in \mathbb{Z}^2} \iota(\sigma(x)) N(x,\sigma)
$$
where
$$
N(x,\sigma) = \sum_{|x-y|=1} \mathbf{1} \{ \sigma(x) \neq \sigma(y) \}
$$
and $\iota(-1)>\iota(0)=\iota(1)>0$ (this last condition will be explained later). This means that at each site a spin is allowed to change at rate $1/2$, independently of any other site, if and only if it does not increase the energy of the system. 

\begin{figure}[htbp]
\includegraphics[viewport=70 600 500 800,clip]{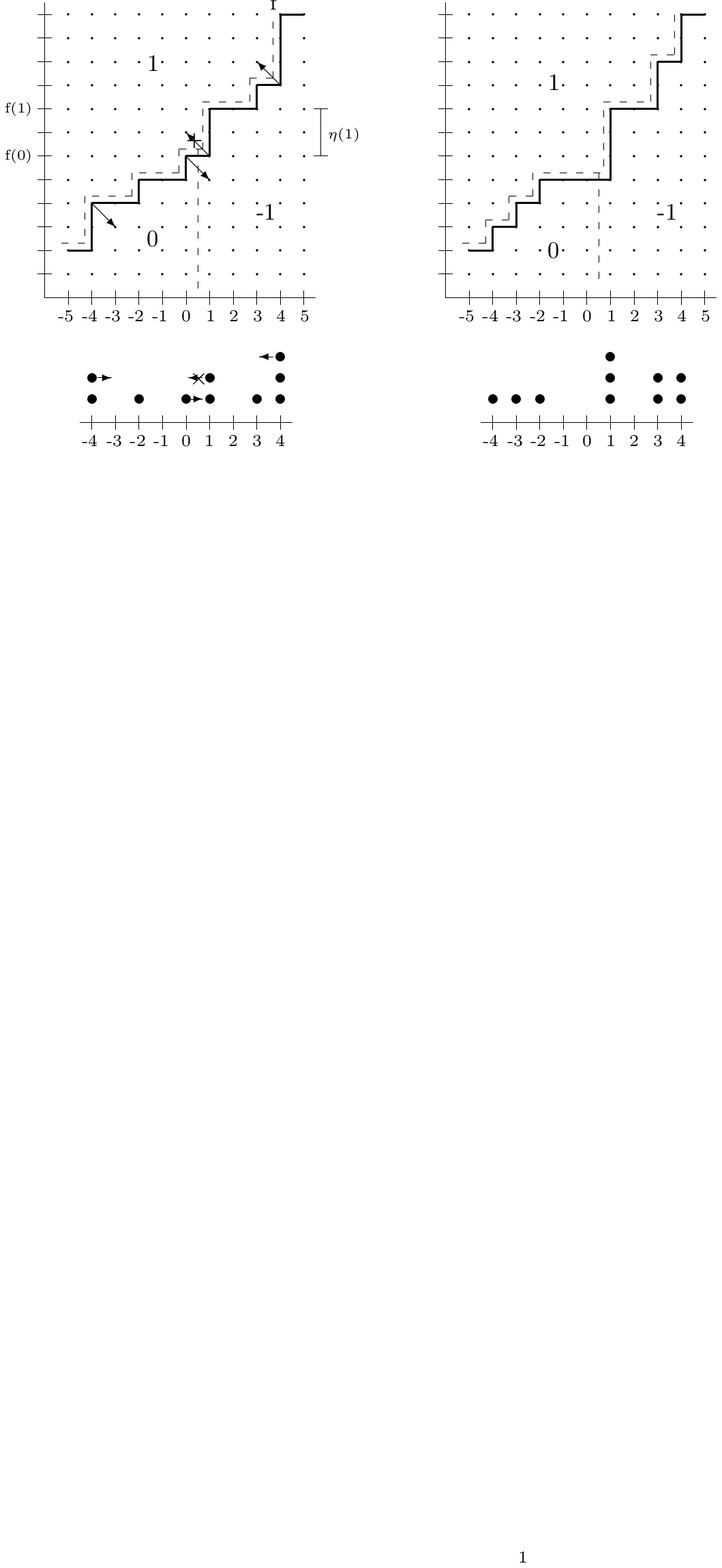}
\caption{\small{On the left: a typical surface configuration and its associated increments configuration with arrows indicating possible transitions and crossed arrows indicating transitions that are not allowed. On the right: modifications to the Potts and zero-range configurations in the left picture following the arrows indications.}}
\label{pic:potts}
\end{figure}

Denote by $\mathcal{I}=\mathcal{I}(\mathbb{Z})$ the collection of non-decreasing functions on $\mathbb{Z}$ and by $\mathcal{A}$ the set of configurations $\sigma$ for which there exists a function $f=f_{\sigma}$ in $\mathcal{I}$ such that, for every $x=(x_1,x_2) \in \mathbb{Z}^2$,
\begin{itemize}
\item[(i)] $\sigma(x)=-1$ if $x_1 >0$, $x_2\le f(x_1)$,
\item[(ii)] $\sigma(x)=0$ if $x_1\le 0$, $x_2\le f(x_1)$,
\item[(iii)] $\sigma(x)=1$ if $x_2>f(x_1)$.
\end{itemize}
This gives a one to one correspondence between $\mathcal{I}$ and $\mathcal{A}$, see figure \ref{pic:potts}. Furthermore, $\mathcal{A}$ is stable under the dynamics induced by the generator $\mathcal{L}_{sp}$. Indeed, fixed a configuration $\sigma \in \mathcal{A}$, by a direct computation of $\mathcal{H}(\sigma^{x,j}) - \mathcal{H}(\sigma)$, it is easy to verify that a jump from $\sigma$ to $\sigma^{x,j}$ is allowed only in one of the following three cases:
\begin{itemize}
\item[(i)] $j=1$ and $x=(x_1,f_{\sigma}(x_1))$ with $f_{\sigma}(x_1)>f_{\sigma}(x_1-1)$,
\item[(ii)] $j=0$ and $x=(x_1,f_{\sigma}(x_1)+1)$ with $x_1<0$ and  $f_{\sigma}(x_1)<f_{\sigma}(x_1+1)$,
\item[(iii)] $j=-1$ and $x=(x_1,f_{\sigma}(x_1)+1)$ with $x_1>0$ and  $f_{\sigma}(x_1)<f_{\sigma}(x_1+1)$.
\end{itemize}
Note that the condition on $\iota$ in the definition of $\mathcal{H}$ is necessary to prevent a jump from $\sigma$ to $\sigma^{(0,f_{\sigma}(0)),-1}$ when $f_{\sigma}(0) > f_{\sigma}(-1)$. Hence we may investigate the evolution of the Markov process $(\sigma_t)$ induced by $\mathcal{L}_{sp}$ and starting from a configuration $\sigma$ in $\mathcal{A}$ through the process $f_t=f_{\sigma_t}$. 

In order to establish the hydrodynamical behavior of the described system, we are now going to introduce some notation and impose some restrictions on the initial conditions of the system. Denote by $\mathbb{N}$, $\mathbb{Z}_-$ and $\mathbb{Z}_+^*$ respectively the sets of non-negative, non-positive and positive integers. It is also clear that a function $f$ in $\mathcal{I}$ is determined by its value at a given site $y \in \mathbb{Z}$ and by its increments, $\{f(x)-f(x-1): \, x \in \mathbb{Z}\}$. Based on this we consider another identification which associates to $f \in \mathcal{I}$ its increments configuration $\eta=\eta_f$ in $\mathbb{N}^{\mathbb{Z}}$ given by
$$
\eta(x)=f(x)-f(x-1), \quad \textrm{for all } x \in \mathbb{Z},
$$
see figure \ref{pic:potts}.

For each $\alpha \ge 0$, denote by $\tilde{\nu}_{\alpha}$ the product measure on $\mathbb{N}^{\mathbb{Z}}$ whose marginals are given by
\begin{equation}
\label{eq:eqmzr}
\tilde{\nu}_{\alpha} \{\eta: \eta(x)=k\} = \frac{1}{1+\alpha} \left( \frac{\alpha}{1+\alpha} \right)^k.
\end{equation}
For each probability measure $m$ on $\mathcal{I}$, denote by $\mathcal{T} m$ the probability measure on $\mathbb{N}^{\mathbb{Z}}$ that corresponds to the distribution of the increments of functions in $\mathcal{I}$. Let $\{m^N: N\ge 1\}$ be a sequence of probability measures on $\mathcal{I}$ satisfying 
\begin{itemize}
\item[(P1)] For every $N\ge1$, $m^N\{f:f(0)=0\}=1$. 
\item[(P2)] Each $\mathcal{T} m^N$ is equal to $(\mathcal{T} m^N)^- \times (\mathcal{T} m^N)^+$, where $(\mathcal{T} m^N)^-$ and $(\mathcal{T} m^N)^+$  are its marginals on $\mathbb{N}^{\mathbb{Z}_{-}}$ and   $\mathbb{N}^{\mathbb{Z}_{+}^*}$ respectively. 
\item[(P3)] The sequence $(\mathcal{T} m^N)$ have marginals on $\mathbb{N}^{\mathbb{Z}_{-}}$ bounded above (resp. below) by $\tilde{\nu}_{\alpha}$ (resp. $\tilde{\nu}_{\lambda}$) for some $0<\lambda <\alpha <\oo$. 
\item[(P4)] There exists a increasing smooth function $\lambda_0: \mathbb{R} \ra \mathbb{R}$ with bounded derivative, such that for each continuous function $G:\mathbb{R} \ra \mathbb{R}$ with compact support and each $\delta>0$,
$$
\lim_{N \ra \oo} m^N \left[ \left| N^{-1} \sum_{x\in \mathbb{Z}} G(x/N) N^{-1}f(x) - \int_{\mathbb{R}} du G(u)\lambda_0(u) \right| \ge \delta \right] = 0 \, .
$$
\end{itemize}

Here (P1) allows us to use the correspondence between $\mathcal{I}$ and $\mathbb{N}^{\mathbb{Z}}$ to study the Potts model through the increments process; (P2) says that the initial condition on sites at left of the origin is independent of the initial condition on sites at the right of the origin allowing us to consider the evolution of the system at left of the origin independently, as we shall discuss later; (P3) implies that the left system will be stochastically dominated by a system with reflection at the origin, which is useful to estimate the density of particles over macroscopic boxes; and (P4) is the usual local equilibrium condition in its weakest form for the sequence of initial probability measures (see \cite{kipnislandim}). 

\medskip 

Let $D(\mathbb{R}_+,\mathcal{I})$ denote the space of right continuous functions with left limits on $\mathcal{I}$ endowed with the skorohod topology. For each probability measure $m$ on $\mathcal{I}$, denote by $\mathbb{P}_m^{sp,N}$ the probability measure on $D(\mathbb{R}_+,\mathcal{I})$ induced by the Markov process $f_t=f_{\sigma_t}$ with generator $\mathcal{L}_{sp}$ speeded up by $N^2$ and initial measure $m$. Our main result is the following:

\begin{theorem} 
\label{theorem:potts}
Fix a sequence of initial measures $\{m^N: N\ge 1\}$ satisfying assumptions (P1)-(P4). For every $\delta>0$
$$
\lim_{N\ra {\oo}} \mathbb{P}_{m^{\scriptscriptstyle{N}}}^{sp,N} \left[ \left| N^{-1} f_t(0) + v_t \right| >\delta \right]=0,
$$
where $v_t$ is given by 
\begin{equation}
\label{eq:v}
v_t = \int_{-\oo}^0 \{\rho_0(u) - \rho(t,u)\} du\, 
\end{equation}
and $\rho$ is the unique weak solution of the nonlinear parabolic equation on $\mathbb{R}_+ \times \mathbb{R}_-$
\begin{equation}
\label{eq:pdedzr}
\left\{ \begin{array}{l}
\partial_t \rho (t,u) =\frac{1}{2} \Delta \Phi(\rho(t,u)), \ \, (t,u) \in (0,+\oo)\times (-\oo,0), \\
\rho(t,0-) = 0, \ \, t \in (0,+\oo),  \\
\rho(0,u) = \partial_u \lambda_0(u), \ \, u \in (-\oo,0),  \\
\end{array} \right. 
\end{equation}
with $\Phi(\rho) = \rho /(1+\rho)$. Moreover, for any continuous function $G:\mathbb{R} \ra \mathbb{R}$ with compact support and any $\delta>0$
$$
\lim_{N\ra {\oo}} \mathbb{P}_{m^{\scriptscriptstyle{N}}}^{sp,N} \left[ \left| \frac{1}{N} \sum_x G(x/N) N^{-1}\{f_t(x)-f_t(0)\} - \int du G(u)\lambda(t,u) \right| \ge \delta \right] =0 \, ,
$$
where $\lambda$ is the unique weak solution of the nonlinear equation
\begin{equation}
\label{eq:pdepm}
\left\{ \begin{array}{l}
\partial_t \lambda(t,u) =\frac{1}{2} \partial_u \Phi(\partial_u \lambda(t,u)), \ \, (t,u) \in (0,+\oo) \times \mathbb{R} \! - \! \{0\}, \\
\partial_u \lambda(t,0-) = 0, \ \, t \in (0,+\oo), \\
\partial_u \Phi(\partial_u \lambda(t,0-)) =\partial_u \Phi(\partial_u \lambda(t,0+)), \ \, t \in (0,+\oo), \\
\lambda(0,u) = \lambda_0(u), \ \, u \in \mathbb{R}. \\
\end{array} \right. 
\end{equation}
\end{theorem} 

\medskip

In the statement above, the integral in (\ref{eq:v}) is to be understood as the limit as $n\ra \oo$ of
$$
\int_{-\oo}^{0} \, du H_n(u) \{ \rho_0(u) - \rho(t,u)\} \, ,
$$
where $H_n(u)=(1+u/n)_+$. Moreover, the precise definitions of weak solutions
for equations (\ref{eq:pdedzr}) and (\ref{eq:pdepm}) will be given later, on sections \ref{sec:weaksolutions} and \ref{sec:exctopotts}, respectively.

Concerning the proof of Theorem \ref{theorem:potts}, as pointed out by Landim,Olla and Volchan \cite{landimollavolchan} the process of the increments of $f_t$ evolves as a zero range process on $\mathbb{Z}$ with asymmetry at the origin. For this process the evolution of the particles on $\mathbb{Z}$ is described by a nearest neighbor, symmetric, space homogeneous zero range process except that a particle is allowed to jump from $0$ to $1$ but not from $1$ to $0$. We should think of this as two coupled process: its restriction to $\mathbb{Z}_-$, the dissipative system; and its restriction to $\mathbb{Z}_+^*$, the absorbing system, which acts as a infinite reservoir for the former.

For the dissipative system, whose evolution is independent of the absorbing system, the hydrodynamic behavior under diffusive scaling was established in \cite{landimollavolchan} with hydrodynamical equation given by the non-linear parabolic equation (\ref{eq:pdedzr}). They also studied the behavior of the total number of particles which leave the system before a fixed time $t>0$. Denoting this number by $X_t$ they proved that $\epsilon X_{\epsilon^{-2} t}$, as $\epsilon \ra 0$, converges in probability to $v_t$ defined in (\ref{eq:v}).

Since the rate at which particles leave the dissipative system is equal to the rate at which particles enter the absorbing system, it is then expected for the coupled process a hydrodynamical behavior under diffusive scaling with hydrodynamic equation:
\begin{equation}
\label{eq:pdezr}
\left\{ \begin{array}{l}
\partial_t \rho(t,u) =\frac{1}{2} \Delta \Phi(\rho(t,u)), \ \, (t,u) \in (0,+\oo) \times \mathbb{R} \! - \! \{0\}, \\
\rho(t,0-) = 0, \ \, t \in (0,+\oo), \\
\partial_u \Phi(\rho(t,0-)) = \partial_u \Phi(\rho(t,0+)), \ \, t \in (0,+\oo), \\
\rho(0,u) = \rho_0(u), \ \, u \in \mathbb{R}. \\
\end{array} \right. 
\end{equation}

However the possible accumulation of particles in a neighborhood of the origin appears as a problem in using a direct approach to establish the hydrodynamical behavior of the absorbing system. To avoid this, we are going to consider a clever transform which maps the absorbing system onto a simple one dimensional nearest-neighbor exclusion process, see Kipnis \cite{kipnis}. Then, we prove the hydrodynamical behavior of this associated exclusion process. This identification is given by an application that associates to each $\eta \in \mathbb{N}^\mathbb{N}$  the configuration $\xi \in \{0,1\}^{\mathbb{Z}_+^*}$ given by
$$
\xi(x) = \left\{
\begin{array}{ll}
1 & , \ \textrm{if} \ x = \sum_{z=1}^n \eta(z) + n \\
0 & , \ \textrm{otherwise}.
\end{array}
\right.
$$
Therefore if $\xi$ is obtained from $\eta$, then $\eta(n)$ is the number of empty sites between (n-1)th and nth particle of $\xi$ for $n > 1$ and the number of empty sites before the first particle of $\xi$ for $n=1$. 

\begin{figure}[htbp]
\includegraphics[viewport=100 620 320 750,clip]{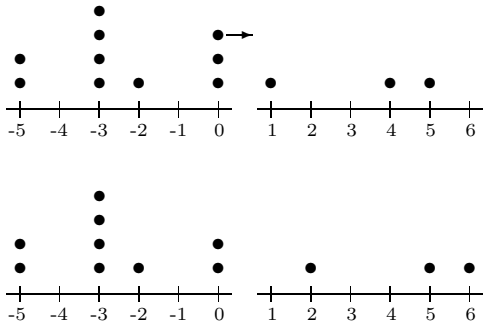}
\caption{\small{A typical configuration in the coupled process followed below it by the configuration obtained when a particle leaves the dissipative system.}}
\label{pic:z-e}
\end{figure}

In this way the absorbing system is mapped onto a process in which each particle jumps as in the nearest neighbor symmetric exclusion process on $\mathbb{Z}_+^*$ with reflection at the origin and superposed to this dynamics, when a particle leaves the dissipative system, the whole system is translated to the right and a new empty site is created at the origin. To a better understand the dynamics that triggers the translation of the system see figure \ref{pic:z-e}. Denoting by 
\begin{equation}
\label{eq:a}
a_t=\frac{dv_t}{dt}=\partial_u \Phi(\rho(t,0-)), \ \, t>0,
\end{equation}
the macroscopic rate at which mass is transfered from the dissipative to the absorbing system, the hydrodynamic equation associated to the coupled exclusion can be derived from an appropriated transformation of the macroscopic profile of the absorbing system, see section \ref{sec:exctopotts}, and is given by
\begin{equation}
\label{eq:pdeep}
\left\{ \begin{array}{l}
\partial_t \zeta(t,u) = \frac{1}{2} \Delta \zeta(t,u) - a_t \partial_u \zeta(t,u), \ \, (t,u) \in (0,+\oo) \times (0,+\oo) \\
\frac{1}{2} \partial_u \zeta(t,0+) = a_t \zeta(t,0+), \ \, t \in (0,+\oo) \\
\zeta(0,u) = \zeta_0(u), \ \, u \in \mathbb{R}_+ 
\end{array} \right. . 
\end{equation}

\smallskip

\begin{remark}From a theoretical point of view, the study of the hydrodynamic behavior for such system is important because of the absence of non-trivial equilibrium measures which does not allow a direct application of the usual methods of proof, see \cite{kipnislandim}. In such cases, the proof becomes particular to each model and few cases have been considered untill now, we refer here to the papers of Chayes-Swindle \cite{chayesswindle}, Landim, Olla and Volchan \cite{landimollavolchan}, Landim and Valle \cite{landimvalle}.  
\end{remark}

\begin{remark}
The coupled process could also be used to describe the behavior of a totally asymmetric tagged particle in a simple symmetric exclusion process under diffusive scaling. This connection have been well exploited in Landim, Olla, Volchan papers \cite{landimollavolchan,landimollavolchan2} and we recommend it to the interested reader.
\end{remark}

\smallskip

The paper has the following structure: In section \ref{sec:hep} we shall consider the coupling between the exclusion process and the dissipative system stating the hydrodynamic limit of the former and section \ref{sec:proof} is devoted to its proof. Finally, in section \ref{sec:exctopotts} we prove Theorem \ref{theorem:potts}. 

\bigskip
\section{Hydrodynamics of the coupled exclusion process}
\label{sec:hep}
\setcounter{equation}{0}

This section is divided in four sub-sections. At first, we present the formal description of the system. After that, we introduce some terminology on weak solutions of the parabolic equations (\ref{eq:pdedzr}) and (\ref{eq:pdeep}) necessary to state the hydrodynamical behavior of the system. In the third part we consider the hypothesis we need on the initial configurations of the system and in the forth we state of the hydrodynamical behavior.

\medskip
\subsection{The system}
\label{sec:system}

The coupled system described informally in the end of the previous section is a Feller process with configuration space $\Omega =\mathbb{N}^{\mathbb{Z}_{-}} \times \{0,1\}^{\mathbb{Z}_{+}^*}$. Denoting by $(\eta,\xi)$ a configuration in $\Omega$, its generator $\mathcal{L}$ may be written as
$$
\mathcal{L} = L + L_b + \widetilde{L}.
$$ 
Here $L$ is related to the motion of particles in the exclusion process: 
$$
L = \sum_{x \ge 1} \{ L_{x,x+1} + L_{x+1,x} \}
$$
where, for every local function $F:\Omega \ra \mathbb{R}$ and every $x,y \ge 1$,
$$
L_{x,y} F(\eta,\xi) = \frac{1}{2} \xi(x) [1-\xi(y)] \ [F(\eta,\xi^{x,y}) - F(\eta,\xi)],
$$
and $\xi^{x,y}$ is the configuration with spins at $x$, $y$ interchanged:
$$
\xi^{x,y}(z)= \left\{ \begin{array}{l}
\xi(y) \, , \ \textrm{if } z=x \\
\xi(x) \, , \ \textrm{if } z=y \\
\xi(z) \, , \ {\textrm{otherwise}} \, . 
\end{array} \right. 
$$
The operator $\widetilde{L}$ is related to the motion of particles in the dissipative system:
$$
\widetilde{L} = \sum_{x \le -1} \{ \widetilde{L}_{x,x+1} + \widetilde{L}_{x+1,x} \}
$$
where, for every local function $F:\Omega \ra \mathbb{R}$ and every $x,y \le 0$,
$$
\widetilde{L}_{x,y}F(\eta,\xi) = \frac{1}{2} g(\eta(x)) [F(\sigma^{x,y} \eta,\xi) - F(\eta,\xi)], 
$$
with $g(k)=\mathbf{1}\{k>0\}$ and
$$
(\sigma^{x,y}\eta)(z)= \left\{ \begin{array}{cl}
\eta(x)-1 & \!\!\!\!\! , \ \textrm{if } z=x \\
\eta(y)+1 & \!\!\!\!\! , \ \textrm{if } z=y \\
\eta(z) & \!\!\!\!\! , \ {\textrm{otherwise}} \, . 
\end{array} \right.  
$$
Finally, $L_b$ is the part of the generator related to the coupling between the systems and thus it describes the jump of a particle at the origin out of the dissipative system and the triggered translation of the whole exclusion process to the right: For every local function $F:\Omega \ra \mathbb{R}$   
$$
L_b F(\eta,\xi) = g(\eta(0)) [F(\eta - \varrho_0,\tau \xi) - F(\eta,\xi)] ,
$$
where $\varrho_x$ stands for the configuration with no particles but one at $x$ and  
$$
(\tau\xi)(x) = \left\{ \begin{array}{cl}
\xi(x-1) & \!\!\!\!\! , \ \textrm{if } x>1 \\
0 & \!\!\!\!\! , \ \textrm{if } x=1 \, .
\end{array} \right. 
$$

\medskip
\subsection{Weak solutions of the hydrodynamical equation}
\label{sec:weaksolutions}

Fix a bounded function $\rho_0:\mathbb{R}_- \ra \mathbb{R}$. A bounded measurable function $\rho:[0,T) \times \mathbb{R}_- \ra \mathbb{R}$ is said to be a weak solution of 
\begin{equation}
\label{eq:pdedzr1}
\left\{ \begin{array}{l}
\partial_t \rho(t,u) = \frac{1}{2} \Delta \Phi(\rho(t,u)) , \ (t,u) \in (0,T) \times (-\oo,0) \\
\rho(t,0-) = 0, \ t \in (0,T) \\
\rho(0,u) = \rho_0(u) , \ u \in \mathbb{R}_- \, ,\\
\end{array} \right. 
\end{equation}
if the following conditions hold:
\begin{itemize}
\item[(a)] $\Phi(\rho(t,u))$ is absolutely continuous in the space variable and $\partial_u\Phi(\rho(t,u))$ is locally square integrable on $(0,T) \times \mathbb{R}_-$ satisfying
$$
\int_0^t \, ds \int_{\mathbb{R}_-} \, du \ e^{u} \{\partial_u\Phi(\rho(s,u))\}^2 < \oo,
\quad \textrm{for every } t>0,
$$
and for every smooth function with compact support $G:[0,T] \times \mathbb{R}_- \ra \mathbb{R}$ vanishing at the origin and for all $0\le t \le T$
$$
\qquad \int_0^t ds \int_{\mathbb{R}_-} du \, G(s,u) \partial_u\Phi(\rho(s,u)) = - \int_0^t ds \int_{\mathbb{R}_-} du \, \partial_uG(s,u) \Phi(\rho(s,u)) \, . \nn
$$
\item[(b)] $\rho(t,0) = 0$ for almost every $0 \le t < T$.
\item[(c)] For every smooth function with compact support $G:\mathbb{R}_- \ra \mathbb{R}$ vanishing at the origin and every $t>0$,
$$
\qquad \int_{\mathbb{R}_-} du \, G(u)\rho(t,u) - \int_{\mathbb{R}_-} du \, G(u)\rho_0(u) =
-\frac{1}{2} \int_0^t ds \int_{\mathbb{R}_-} du \, G^{\prime}(u) \partial_u\Phi(\rho(s,u)). 
$$  
\end{itemize}

\smallskip

\no For a uniqueness result for equation (\ref{eq:pdedzr1}) see \cite{landimollavolchan}.

\medskip

Now, for a fixed bounded function $\zeta_0:\mathbb{R}_+ \ra \mathbb{R}$. A bounded measurable function $\zeta:[0,T) \times \mathbb{R}_+ \ra \mathbb{R}$ is said to be a weak solution of 
\begin{equation}
\label{eq:pdeep1}
\left\{ \begin{array}{l}
\partial_t \zeta(t,u) = \frac{1}{2} \Delta \zeta(t,u) - a_t \partial_u \zeta(t,u), \ \, (t,u) \in (0,T) \times (0,+\oo) \\
\frac{1}{2} \partial_u \zeta(t,0+) = a_t \zeta(t,0+), \ \, t \in (0,T) \\
\zeta(0,u) = \zeta_0(u), \ \, u \in \mathbb{R}_+ 
\end{array} \right. ,
\end{equation}
where $a:(0,T) \ra \mathbb{R}_+$ is a bounded measurable function, if
\begin{itemize}
\item[(a)] $\zeta(t,u)$ is absolutely continuous in the space variable and $\partial_u\zeta(t,u)$ is a locally square integrable function on $(0,T) \times \mathbb{R}_+$ such that for all $0\le t \le T$ and for every smooth function $G:[0,T] \times \mathbb{R}_+ \ra \mathbb{R}$ with compact support
\begin{eqnarray}
\lefteqn{ \qquad \quad \int_0^T ds \int_{\mathbb{R}_+} du \, G(s,u) \partial_u\zeta(s,u) = } \nn \\
& & \qquad \quad = - \int_0^T ds \int_{\mathbb{R}_+} du \, \partial_uG(s,u) \zeta(s,u) - \lim_{\epsilon \ra 0} \int_0^T ds \, G(s,0) \frac{1}{\epsilon} \int_0^\epsilon \zeta(s,u) du \, . \nn
\end{eqnarray}
\item[(b)] For every smooth function with compact support $G:[0,T] \times \mathbb{R}_+ \ra \mathbb{R}$ and every $t\in [0,T]$,
\begin{eqnarray}
\lefteqn{\qquad \int_{\mathbb{R}_+} du \, G(t,u)\zeta(t,u) - \int_{\mathbb{R}_+} du \, G(0,u)\zeta_0(u) = \int_0^t ds \int_{\mathbb{R}_+} du \, \partial_s G(s,u) \zeta(s,u) +
} \nn \\
& & \qquad \qquad + \ \int_0^t ds \left\{  
-\frac{1}{2} \int_{\mathbb{R}_+} du \, G^{\prime}(s,u) \partial_u\zeta(s,u) +a_s \int_{\mathbb{R}_+} du \, G^{\prime}(s,u) \zeta(s,u) \right\} . \nn
\end{eqnarray}  
\end{itemize}

\begin{remark}The method we apply to prove the hydrodynamical behavior of the coupled exclusion process, see section \ref{sec:initialmeasures} below, requires uniqueness of solutions of equation (\ref{eq:pdeep1}). On the other hand existence and regularity of solutions of equation (\ref{eq:pdeep1}) follows from the proof of the hydrodynamical behavior itself, which is described in Theorem \ref{theorem:hbep} on section \ref{sec:hbep}, this should be clear later on, see also \cite{kipnislandim}. The required uniqueness for the equation is given in Theorem 5.1 of chapter 3 in \cite{lsu} together with the remark at end of section 5 in the same chapter of the book.
\end{remark}

\medskip
\subsection{Hypothesis on the initial measures}
\label{sec:initialmeasures}

Given any two measures $\mu$, $\nu$ on $\Omega$, we denote by $H(\mu|\nu)$ the relative entropy of $\mu$ with respect to $\nu$:
$$
H(\mu|\nu)=\sup_f \left\{ \int fd\mu - \log \int e^f d\nu \right\} ,
$$
where the supremum is carried over all bounded continuous real functions on $\Omega$. From this variational formula, we have the so called entropy inequality, i.e., for every bounded continuous real function $f$ on $\Omega$
\begin{equation}
\label{eq:entinq}
\int fd\mu \le \frac{H(\mu|\nu)}{a} + \frac{1}{a} \log \int e^{af} d\nu
\end{equation}
for every constant $a>0$. Recall also that if $\mu$ is absolutely continuous with respect to $\nu$, then
$$
H(\mu|\nu) = \int \log \frac{d\mu}{d\nu} d\mu \, ,
$$
see for instance appendix 1.8 in \cite{kipnislandim}.

For a measure $\mu$ on $\Omega$, denote by $\mu^-$ and $\mu^+$ its marginals on $\mathbb{N}^{\mathbb{Z}_{-}}$ and $\{0,1\}^{\mathbb{Z}_+^*}$, respectively. Let $\mathcal{P}_{\pm}(\Omega)$ be the space of probability measures $\mu$ on $\Omega$ that can be written as $\mu=\mu^+ \times \mu^-$. For $0<\alpha<1$, let $\nu_\alpha$ denote the Bernoulli product measure of parameter $\alpha$ on $\{0,1\}^{\mathbb{Z}_+^*}$. 

Fix a sequence of probability measures $\{ \mu^{N} : N\ge 1 \}$ on $\mathcal{P}_{\pm}(\Omega)$. To prove the hydrodynamical behavior of the system we shall assume that
\begin{itemize}
\item[(E1)] The sequence $(\mu^{N,-})$ is bounded above (resp. below) by $\tilde{\nu}_{\lambda_1}$ (resp. $\tilde{\nu}_{\lambda_2}$) for some $0<\lambda_1<\lambda_2$.
\item[(E2)] There exists a $\beta>0$ such that $H(\mu^{N,-}|\tilde{\nu}_\beta) \le CN$ for some constant $C>0$, where $\tilde{\nu}_\beta$ is defined in (\ref{eq:eqmzr}).
\item[(E3)] The sequence $\{\mu^{N,-}, N\ge 1\}$ is associated to a bounded initial profile $\rho_0: \mathbb{R_-} \ra \mathbb{R}$, i.e., for each $\delta>0$ and each continuous function $G:\mathbb{R}_- \ra \mathbb{R}$ with compact support
$$
\lim_{N\ra {\oo}} \mu^N \left[ \left| \frac{1}{N} \sum_{x\le 0} G(x/N)\eta(x) - \int_{-\oo}^0 duG(u)\rho_0(u) \right| \ge \delta \right] = 0.
$$ 
\item[(E4)]  The sequence $\{\mu^{N,+}, N\ge 1\}$ is associated to a bounded initial profile $\zeta_0: \mathbb{R_+} \ra \mathbb{R}$, i.e., for each $\delta>0$ and each continuous function $G:\mathbb{R}_+ \ra \mathbb{R}$ with compact support
$$
\lim_{N\ra {\oo}} \mu^N \left[ \left| \frac{1}{N} \sum_{x\ge 1} G(x/N)\xi(x) - \int_0^{+\oo} duG(u)\zeta_0(u) \right| \ge \delta \right] = 0.
$$ 
\end{itemize}
\smallskip
The first three assumptions are used by Landim, Olla and Volchan in \cite{landimollavolchan} to establish the hydrodynamical behavior of the dissipative system. The condition (E4) is the usual law of large numbers imposed on the empirical measure at time $0$ for the coupled exclusion.

We include one last condition in the previous list:
\begin{itemize}
\item[(E5)] There exists $0< \alpha <1$ such that $H(\mu^{N,+}|\nu_\alpha) \le CN$ for some constant $C>0$.  
\end{itemize}
This condition will be required in section \ref{sec:proof} to prove Lemma \ref{lemma:replacementbymean}. However, the necessity of condition (E5) is removed is section \ref{sec:rementcond}. 

The proof of the hydrodynamical behavior of the dissipative system in \cite{landimollavolchan} is an adaptation of the entropy method introduced by 
Guo, Papanicolaou and Varadhan in \cite{guopapanicolaouvaradhan}, see also \cite{kipnislandim}. Such an adaptation requires the introduction of appropriate reference measures that plays the role of the missing equilibrium measures. For the dissipative system the reference measures are the product measures $\tilde{\nu}_{\gamma(\cdot)}^N$, $N \ge 1$, on $\mathbb{N}^{\mathbb{Z}_-}$ with marginals given by
$$
\tilde{\nu}_{\gamma(\cdot)}^N \{\eta : \eta(x) = k\} = \tilde{\nu}_{\Phi^{-1}(\gamma_x^N)} \{\eta : \eta(x) = k\} = (1-\gamma_x^N) (\gamma_x^N)^k \, .
$$
where for each $x \le 0$ 
$$
\gamma_x^N := 
\left\{ \begin{array}{cl}
\frac{\beta (1+x)}{N} & \!\!\!\!\! , \ \textrm{for } -N+1 \le x \le 0 \\
\beta & \!\!\!\!\! , \ \textrm{for } x \le -N .
\end{array} \right.
$$
with $\beta >0$ satisfying (E2). For these measures the entropy bound is preserved , which means (see \cite{landimollavolchan}) that
\begin{equation}
\label{eq:entbound}
H(\mu^{N,-} | \tilde{\nu}_{\gamma(\cdot)}^N) \le C N \, .
\end{equation}

Although we do not need an entropy production estimate for the exclusion process we shall adapt some of the estimates found in \cite{landimollavolchan} that requires adequate reference measures for this system. Taking the reference measure for the absorbing system as the canonical measure $\nu_\alpha$, the reference measures for the whole system are taken as the product of the reference measures for the dissipative system with $\nu_\alpha$, i.e., $\nu_{\alpha,\gamma}^N := \tilde{\nu}_{\gamma(\cdot)}^N \times \nu_\alpha$.

\medskip
\subsection{The hydrodynamical behavior}
\label{sec:hbep}
Let $D(\mathbb{R}_+,\Omega)$ denote the space of right continuous functions with left limits on $\Omega$ endowed with the skorohod topology. For each probability measure $\mu$ on $\Omega$, denote by $\mathbb{P}_{\mu}^N$ the probability measure on $D(\mathbb{R}_+,\Omega)$ induced by the Markov process $(\eta_t,\xi_t)$ with generator $\mathcal{L}$ speeded up by $N^2$ and with initial measure $\mu$. The hydrodynamical behavior of the exclusion process is given by the following result:

\begin{theorem} 
\label{theorem:hbep}
Fix a sequence of initial measures $\{\mu^N, N>1\}$ on $\mathcal{P}_{\pm}(\Omega)$ satisfying (E1)-(E4) with strictly positive initial profile $\zeta_0: \mathbb{R_+} \ra \mathbb{R}$ bounded above by $1$. Then, for any continuous $G: \mathbb{R_+} \ra \mathbb{R}$ with compact support, any $\delta>0$ and $0<t<T$ 
$$
\lim_{N\ra {\oo}} \mathbb{P}_{\mu^{\scriptscriptstyle{N}}}^N \left[ \left| \frac{1}{N} \sum_{x\ge 1} G(x/N)\xi_t(x) - \int duG(u)\zeta(t,u) \right| \ge \delta \right] = 0
$$
where $\zeta$ is the unique solution of (\ref{eq:pdeep1}), with $a_t=\partial_u \Phi(\rho(t,0-))$ for $\rho$ being the unique solution of (\ref{eq:pdedzr1}).
\end{theorem} 

\medskip
\section{The proof of the Hydrodynamic limit}
\label{sec:proof}
\setcounter{equation}{0}

Denote by $\mathcal{M}=\mathcal{M}(\mathbb{R})$, the space of positive Radon measures on $\mathbb{R}$ endowed with the vague topology. Integration of a function $G$ with respect to a measure $\pi$ in $\mathcal{M}$ will be denoted $\< \pi,G\rangle$. To each configuration $(\eta,\xi) \in \Omega$ and each $N\ge 1$ we associate the empirical measure $\pi^N + \tilde{\pi}^N$ in $\mathcal{M}$, where
$$
\pi^N = \frac{1}{N} \sum_{x\ge 1} \xi(x) \delta_{x/N}  \qquad \textrm{and} \qquad \tilde{\pi}^N = \frac{1}{N} \sum_{x\le 0} \eta(x) \delta_{x/N} \, .
$$
Let $D([0,T],\mathcal{M})$ denote the space of right continuous functions with left limits on $\mathcal{M}$ endowed with the Skorohod topology and on the space of probability measures on $D([0,T],\mathcal{M})$ we also consider the vague topology. For each probability measure $\mu$ on $\Omega$, denote by $\mathbb{Q}_{\mu}^N$ (resp. $\tilde{\mathbb{Q}}_{\mu}^N$) the probability measure on $D([0,T],\mathcal{M})$ induced by $\mathbb{P}_{\mu}^N$ and the empirical measure $\pi^N$ (resp. $\tilde{\pi}^N$).

Theorem \ref{theorem:hbep} states that the sequence $\mathbb{Q}_{\mu^N}^N$ converges weakly, as $N \ra \oo$, to the probability measure concentrated on absolutely continuous trajectories $\pi(t,du)=\zeta(t,u)du$ whose density is the solution of (\ref{eq:pdeep1}) (see \cite{kipnislandim}). The proof consists of showing tightness of $\mathbb{Q}_{\mu^N}^N$, that all of its limit points are concentrated on absolutely continuous paths which are weak solutions of (\ref{eq:pdeep1}) and uniqueness of solutions of this equation.

We have already discussed uniqueness of weak solutions of (\ref{eq:pdeep1}) in section \ref{sec:weaksolutions}.

Note that all limit points of the sequence $\mathbb{Q}_{\mu}^N$ are concentrated on absolutely continuous measures since the total mass on compact intervals of the empirical measure $\pi^N$ is bounded by the size of the interval plus $1/N$.

In order to show that all limit points of the sequence $\mathbb{Q}_{\mu^{\scriptscriptstyle{N}}}^N$ are concentrated on weak solutions of (\ref{eq:pdeep1}) we will need the following result:

\begin{lemma}
\label{lemma:weaksol}
For every smooth function $G:[0,T] \times \mathbb{R}_+ \ra \mathbb{R}$ with compact support and $\delta>0$
\begin{eqnarray}
\lefteqn{ \!\!\!\!\!\!\!\!\!
\limsup_{\epsilon \ra 0} \limsup_{N \ra \oo} \mathbb{Q}_{\mu^{\scriptscriptstyle{N}}}^N \left[ \sup_{0\le t\le T} \left| \<\pi_t^N,G\rangle - \<\pi_0^N,G\rangle - \int_0^t \{ \< \pi_s^N, \partial_s G \rangle + \frac{1}{2} \< \pi_s^N,\Delta G \rangle + \right. \right.} \nn \\ 
& & \qquad \qquad \qquad \left. \left. + \frac{1}{2} \nabla G(s,0) \<\pi_s^N, \textbf{1}[0,\epsilon] \rangle + a_s \< \pi_s^N, \nabla G \rangle  \} ds  \right| > \delta \right] = 0 \, . \nn
\end{eqnarray}
\end{lemma}

\smallskip

We shall divide the proof of Theorem \ref{theorem:hbep} in four parts: We start proving tightness in section \ref{sec:tightness}. The section \ref{sec:weaksol} is devoted to prove of Lemma \ref{lemma:weaksol} under a condition on the entropy of the system with respect to an equilibrium measure for the system with reflexion at the origin. From Lemma \ref{lemma:weaksol}, to conclude the proof that all limit points of the sequence $\mathbb{Q}_{\mu^{\scriptscriptstyle{N}}}^N$ are concentrated on weak solutions of (\ref{eq:pdeep1}), we have to justify an integration by parts to obtain conditions (a) and (b) in the definition of weak solutions of (\ref{eq:pdeep1}). This is consequence of an energy estimate which is the content of section \ref{sec:energyest}. We finish the proof of Theorem \ref{theorem:hbep} removing the imposed condition on the entropy at section \ref{sec:rementcond}.

\subsection{Tightness}
\label{sec:tightness}  

The sequence $\mathbb{Q}_{\mu}^N$ is tight in the space of probability measures on $D([0,T],\mathcal{M})$, if for each smooth function with compact support $G:\mathbb{R}_+ \ra \mathbb{R}$, $\<\pi_t^N,G\rangle$ is tight as a random sequence on $D(\mathbb{R}_+,\mathbb{R})$. Now fix such a function, denote by $\mathcal{F}_t = \sigma ((\tilde{\pi}_s,\pi_s), s\le t)$, $t \ge 0$, the natural filtration on $D([0,T],\mathcal{M})$, and by $\mathcal{T}_T$ the family of stopping times bounded by $T$. According to Aldous \cite{aldous}, to prove tightness for $\<\pi_t^N,G\rangle$ we have to verify the following two conditions:
\begin{itemize}
\item[(i)] The finite dimensional distributions of $\<\pi_t^N,G\rangle$ are tight;
\item[(ii)] for every $\epsilon >0$ 
\begin{equation}
\label{eq:aldous}
\lim_{\gamma \ra 0} \limsup_{N \ra \oo} \sup_{\tau \in \mathcal{T}_T} \sup_{\theta \le \gamma} \,
\mathbb{P}_{\mu^{\scriptscriptstyle{N}}}^N \left[ |\<\pi^N_\tau,G\rangle - \<\pi^N_{\tau + \theta},G\rangle | > \epsilon \right] = 0 \, .
\end{equation}
\end{itemize}

Condition (i) is a trivial consequence of the fact that the empirical measure has finite total mass on any compact interval. In order to prove condition (ii), Let us first consider for each smooth function $H:[0,T] \times \mathbb{R}_+ \ra \mathbb{R}¨$ the associated $(\mathcal{F}_t)$-martingale  vanishing at the origin
\begin{equation}
\label{eq:martingal}
M_t^{H,N} = \<\pi_t^N,H\rangle - \<\pi_0^N,H\rangle - \int_0^t (\partial_s + N^2 (L+L_b)) \<\pi_s^N,H\rangle ds \, .
\end{equation}
In (ii) the function $G$ does not depend on the time, however the martingale is defined for functions varying on time, which we use in the proof of lemma \ref{lemma:weaksol}.  An elementary computation shows that $(\partial_s + N^2(L+L_b)) \<\pi^N,H\rangle$ is given explicitly by 
\begin{equation}
\label{eq:integrand}
\< \pi^N,\partial_s H \rangle +
\frac{1}{2} \< \pi^N,\Delta_N H \rangle + \frac{1}{2} \nabla_N H(s,0) \xi(1) +  N g(\eta(0)) \< \pi^N, \nabla_N H \rangle
\end{equation}
where $\Delta_N$ and $\nabla_N$ denote respectively the discrete Laplacian and gradient:
$$
\Delta_NH(s,x/N) = N^2 \{H(s,(x+1)/N) + H(s,(x-1)/N) - 2H(s,x/N)\} \, ,
$$
$$
\nabla_NH(s,x/n) = N\{H(s,(x+1)/N) - H(s,x/N)\}. 
$$
We also derive an explicit formula for the quadratic variation of $M^{H,N}_t$, see Lemma 5.1 in Appendix 1 of \cite{kipnislandim}. It is given by
\begin{eqnarray}
\label{eq:quavar}
\lefteqn{ \!\!\!\!\!\!\!\!\!\!\!\!\!\!
\<M^{H,N}\rangle_t = \int_0^t ds \frac{1}{N^2} \left\{ \sum_{x,y \ge 1 \atop |x-y|=1} \nabla_NH(s,x \land y/N)^2 \xi_s(x)[1-\xi_s(y)] + \right.}\nn \\
& & \qquad \qquad  \left. + g(\eta_s(0)) \sum_{x,y\ge 1} \nabla_NH(s,x/N) \nabla_NH(s,y/N) \xi_s(x) \xi_s(y) \} \right\} .
\end{eqnarray}

Therefore, for $G$ not depending on the time as before
$$
\<\pi^N_{\tau + \theta},G\rangle - \<\pi^N_\tau,G\rangle = M_{\tau+\theta}^{G,N} - M_{\tau}^{G,N} + \int_\tau^{\tau + \theta} N^2 \mathcal{L}\<\pi_s^N,G\rangle ds \, .
$$
From the previous expression and Chebyshev inequality, (ii) follows from
\begin{equation}
\label{eq:decmart}
\lim_{\gamma \ra 0} \limsup_{N \ra \oo} \sup_{\tau \in \mathcal{T}_T} \sup_{\theta \le \gamma} \,
\mathbb{E}_{\mu^{\scriptscriptstyle{N}}}^N \left[ \left| M_{\tau+\theta}^{G,N} - M_{\tau}^{G,N} \right| \right] = 0 
\end{equation}
and
\begin{equation}
\label{eq:decint}
\lim_{\gamma \ra 0} \limsup_{N \ra \oo} \sup_{\tau \in \mathcal{T}_T} \sup_{\theta \le \gamma} \,
\mathbb{E}_{\mu^{\scriptscriptstyle{N}}}^N \left[ \left| \int_\tau^{\tau + \theta} N^2 \mathcal{L}\<\pi_s^N,G\rangle ds \right| \right] = 0 \, .
\end{equation} 

\smallskip

Now we show (\ref{eq:decmart}) and (\ref{eq:decint}) and complete the proof of tightness.

\smallskip

\no \textbf{Proof of (\ref{eq:decmart}):} From the optional stopping theorem and the martingale property
$$
\mathbb{E}_{\mu^{\scriptscriptstyle{N}}}^N \left[ (M_{\tau+\theta}^{G,N} - M_{\tau}^{G,N})^2 \right] = \mathbb{E}_{\mu^{\scriptscriptstyle{N}}}^N \left[ \<M^{G,N}\rangle_{\tau+\theta} - \<M^{G,N}\rangle_\tau \right] \, .
$$
Hence, applying formula (\ref{eq:quavar}), by the Taylor expansion for G, we have that 
\begin{eqnarray}
\label{eq:aldous2}
\lefteqn{\!\!\!\!\!\!\!\!\!\!\!\!\!\!\!\!\!\!\!\!\!\!\!\!\!
\mathbb{E}_{\mu^{\scriptscriptstyle{N}}}^N \left[ (M_{\tau+\theta}^{G,N} - M_{\tau}^{G,N})^2 \right] \le \frac{C(G)}{N} \left( \theta + \mathbb{E}_{\mu^{\scriptscriptstyle{N}}}^N \left[ \int_\tau^{\tau + \theta} N g(\eta_s(0)) ds \right] \right)} \nn \\
& & \le \frac{C(G)}{N} \left( \theta +  \mathbb{E}_{\mu^{\scriptscriptstyle{N}}}^N \left[ \int_0^{T+\theta} N g(\eta_s(0)) ds \right]  \right).
\end{eqnarray}
where the previous inequality holds because $\tau \in \mathcal{T_T}$.
Recall from the proof of Lemma 3.8 in \cite{landimollavolchan} that for every $0 < t \le T$,
\begin{equation}
\label{eq:basicdecay} 
\mathbb{E}_{\mu^{\scriptscriptstyle{N}}}^N \left[ \int_0^t g(\eta_s(0)) ds \right] \le \frac{t}{\sqrt{N}} .
\end{equation}
Therefore, from (\ref{eq:aldous2}), we have 
$$
\sup_{\tau \in \mathcal{T}_T} \mathbb{E}_{\mu^{\scriptscriptstyle{N}}}^N \left[ (M_{\tau+\theta}^{G,N} - M_{\tau}^{G,N})^2 \right] \le \frac{C(G)(T+\theta)}{\sqrt{N}}
$$ 
and (\ref{eq:decmart}) holds. $\square$

\smallskip

\no \textbf{Proof of (\ref{eq:decint}):} 
From formula (\ref{eq:integrand}) and the Taylor expansion for $G$ we obtain that
\begin{equation}
\label{eq:aldous1}
\mathbb{E}_{\mu^{\scriptscriptstyle{N}}}^N \left[ \left| \int_\tau^{\tau + \theta} \! \! N^2 \mathcal{L}\<\pi_s^N,G\rangle ds \right| \right] \le C(G) \left( \theta + \mathbb{E}_{\mu^{\scriptscriptstyle{N}}}^N \left[ \int_\tau^{\tau + \theta} \! \! N g(\eta_s(0)) ds \right] \right).
\end{equation}
Hence (\ref{eq:decint}) follows from (\ref{eq:aldous1}) if 
\begin{equation}
\label{eq:stdecay} 
\lim_{\theta \ra 0} \limsup_N \, \sup_{\tau \in \mathcal{T}_T} \mathbb{E}_{\mu^{\scriptscriptstyle{N}}}^N \left[ \int_\tau^{\tau + \theta} N g(\eta_s(0)) ds \right] =0 .
\end{equation}
We postpone the proof (\ref{eq:stdecay}), assuming it we have (\ref{eq:decint}). $\square$

\medskip

It remains to prove (\ref{eq:stdecay}) to complete the proof of (\ref{eq:aldous}). However we first show that in fact the expectation in (\ref{eq:basicdecay}) is of order $O(N^{-1})$, i.e.,
\begin{equation}
\label{eq:decay}
\sup_N \mathbb{E}_{\mu^{\scriptscriptstyle{N}}}^N \left[ \int_0^t N g(\eta_s(0)) ds \right] < \oo .
\end{equation}
The reason is that (\ref{eq:decay}) is required in the next sections and its proof is similar to that of  (\ref{eq:stdecay}). 

\smallskip

\no \textbf{Proof of \ref{eq:decay}:} We introduce a second class of martingales. Let $G:\mathbb{R}_{-} \ra \mathbb{R}$ be a smooth function with compact support and denote by $\tilde{M}^{G,N}$ the $(\mathcal{F}_t)$-martingale
$$
\tilde{M}_t^{G,N} = \<\tilde{\pi}_t^N,G\rangle - \<\tilde{\pi}_0^N,G\rangle - \int_0^t N^2 (\tilde{L}+L_b) \<\tilde{\pi}_s^N,G\rangle ds \, .
$$
By a straightforward computation we have that $N^2 (\tilde{L}+L_b) \<\tilde{\pi}^N,G\rangle$ is equal to
$$
\frac{1}{2N} \sum^{-1}_{x=-\oo} \Delta_N G(x/N) g(\eta(x)) - \frac{1}{2} \nabla_N G(-1/N) g(\eta(0)) - NG(0)g(\eta(0)) \, .
$$
Moreover, the quadratic variation $\<\tilde{M}^{G,N}\rangle_t$ is given by
$$
\int_0^t ds \frac{1}{2N^2} \left\{ \sum_{x,y \le 0 \atop |x-y|=1} \nabla_NG(x \land y/N)^2 [g(\eta_s(x))+g(\eta_s(y))] + N^2 G(0)^2 g(\eta_s(0)) \right\} \, .
$$
For $l \in \mathbb{N}$, let $H_l: \mathbb{R}_- \ra \mathbb{R}$ be the function $H_l(u)=(1+u/l)_+$, $u\le 0$. Then 
$$
\tilde{M}^{H_l,N}_t + \<\tilde{\pi}^N_0,H_l\rangle - \<\tilde{\pi}^N_t,H_l\rangle = \int_0^t ds \left[ N g(\eta_s(0)) + \frac{g(\eta_s(0)) - g(\eta_s(-lN))}{2l} \right] .
$$
Since $g\le 1$, for each configuration $\eta \in \mathbb{N}^{\mathbb{Z}_-}$,
\begin{equation}
\label{eq:nghl}
\int_0^t N g(\eta_s(0)) ds = \lim_{l \ra \oo} \left\{ \tilde{M}^{H_l,N}_t + \<\tilde{\pi}^N_0,H_l\rangle - \<\tilde{\pi}^N_t,H_l\rangle \right\} .
\end{equation}
By Fatou's Lemma, we have that the expectation of the right hand side term in the previous equality is bounded above by
$$
\sup_N \limsup_{l \ra \oo} \mathbb{E}_{\mu^{\scriptscriptstyle{N}}}^N \left[ | \tilde{M}^{H_l,N}_t + \<\tilde{\pi}^N_0,H_l\rangle - \<\tilde{\pi}^N_t,H_l\rangle | \right]. 
$$
Note that
\begin{eqnarray}
\lefteqn{
\mathbb{E}_{\mu^{\scriptscriptstyle{N}}}^N \left[ | \tilde{M}_t^{H_l,N} | \right]^2 \ \le \ \mathbb{E}_{\mu^{\scriptscriptstyle{N}}}^N \left[ | \tilde{M}_t^{H_l,N} |^2 \right] \ = \ \mathbb{E}_{\mu^{\scriptscriptstyle{N}}}^N \left[ \< \tilde{M}^{H_l,N} \rangle_t \right] } \nn \\
& \quad \qquad = & \mathbb{E}_{\mu^{\scriptscriptstyle{N}}}^N \left[ \int_0^t ds \left\{ g(\eta_s(0)) + 2 \! \! \sum_{-lN\le x \le -1} \! \! \frac{g(\eta_s(x)) + g(\eta_s(x+1))}{(lN)^2} \right\} \right] \nn \\
& \quad \qquad \le &\mathbb{E}_{\mu^{\scriptscriptstyle{N}}}^N \left[ \int_0^t g(\eta_s(0)) ds \right] + \frac{4t}{lN} . \nn
\end{eqnarray}
In particular, by (\ref{eq:basicdecay}),
\begin{equation}
\label{eq:martdecay}
\lim_{N \ra \oo} \sup_l \mathbb{E}_{\mu^{\scriptscriptstyle{N}}}^N \left[ \left| \tilde{M}_t^{H_l,N} \right| \right] = 0 \, ,
\end{equation}
so that (\ref{eq:decay}) holds if
$$
\sup_N \limsup_l \mathbb{E}_{\mu^{\scriptscriptstyle{N}}}^N \left[  | \<\tilde{\pi}_t^N,H_l\rangle - \<\tilde{\pi}_0^N,H_l\rangle | \right] < \oo .
$$  
To see this, fix $C>0$ and a sequence of continuous functions $\{ G_l:\mathbb{R}_- \ra \mathbb{R}\}$ bounded by one and vanishing at the origin such that, $G_l$ is smooth on $(-l,0)$, $G_l = H_l$ on $(-\oo,-C]$, and $\{\nabla G_l: (C,0) \ra \mathbb{R}\}$, $\{\Delta G_l: (C,0) \ra \mathbb{R}\}$ are uniformly bounded families of functions (It is straightfoward to obtain such functions, we let this to the reader). Then $| \<\tilde{\pi}_t^N,H_l\rangle - \<\tilde{\pi}_0^N,H_l\rangle |$ is bounded above by
$$
| \<\tilde{\pi}_t^N,H_l\rangle - \<\tilde{\pi}_t^N,G_l\rangle | 
+ | \<\tilde{\pi}_t^N,G_l\rangle - \<\tilde{\pi}_0^N,G_l\rangle | +
| \<\tilde{\pi}_0^N,H_l\rangle - \<\tilde{\pi}_0^N,G_l\rangle | .
$$
Using the martingale $\tilde{M}^{G_l,N}$ and its quadratic variation, we verify by usual computations that the expectation of the middle term at the right hand side of this equation is uniformly bounded in both $N$ and $l$. The other two terms are bounded respectively by
$$
\frac{1}{N} \sum_{x=-CN}^{0} \eta_t^N(x) \quad \textrm{and} \quad \frac{1}{N} \sum_{x=-CN}^{0} \eta_0^N(x) ,
$$
whose expectation is also uniformly bounded. To prove this last statement, use the fact that $\mu^{N,-} \le \tilde{\nu}_\alpha$ to construct a coupling, between the dissipative system and the nearest-neighbor, symmetric,  space-homogeneous zero-range process with reflection at the origin, which preserves the stochastic order and recall that for this last system $\tilde{\nu}_\alpha$, defined in (\ref{eq:eqmzr}), is an equilibrium state. 
Therefore (\ref{eq:decay}) holds. $\square$

\smallskip

\no \textbf{Proof of (\ref{eq:stdecay}):} As in (\ref{eq:nghl}),  we have that 
\begin{equation}
\label{eq:nghl2}
\int_{\tau}^{\tau + \theta} N g(\eta_s(0)) ds = \lim_{l \ra \oo} \left\{ \tilde{M}^{H_l,N}_{\tau + \theta} - \tilde{M}^{H_l,N}_\tau + \<\tilde{\pi}^N_\tau ,H_l\rangle - \<\tilde{\pi}^N_{\tau + \theta},H_l\rangle \right\} .
\end{equation}
On the one hand, as in the proof of (\ref{eq:decmart}), obtain by quadratic variation of $\tilde{M}^{H_l,N}$ 
\begin{eqnarray}
\lefteqn{ \lim_{N \ra \oo} \limsup_{l\ra \oo} \sup_{\tau \in \mathcal{T}_T} \, \mathbb{E}_{\mu^{\scriptscriptstyle{N}}}^N \left[ | \tilde{M}^{H_l,N}_{\tau + \theta} - \tilde{M}^{H_l,N}_\tau | \right]^2 \le } \nn \\
& \qquad \qquad \le & \lim_{N \ra \oo} \limsup_{l\ra \oo} \, \left\{ \mathbb{E}_{\mu^{\scriptscriptstyle{N}}}^N \left[ \int_0^T g(\eta_s(0)) ds \right] + \frac{4T}{lN} \right\} = 0 \, . \nn
\end{eqnarray}
On the other hand, $| \<\tilde{\pi}_{\tau + \theta}^N,H_l\rangle - \<\tilde{\pi}_{\tau}^N,H_l\rangle |$
is bounded above by
$$
| \<\tilde{\pi}_{\tau + \theta}^N,G_l\rangle - \<\tilde{\pi}_{\tau}^N,G_l\rangle | + 
\frac{1}{N} \sum_{x=-CN}^{0} \eta_{\tau +\theta}^N(x) + \frac{1}{N} \sum_{x=-CN}^{0} \eta_{\tau}^N(x) \, ,  
$$
with $G_l$ taken as in the proof of (\ref{eq:decay}).
Using again the explicit formulas for $\tilde{M}_{G_l,N}$ and for its quadratic variation, we show that 
$$
\lim_{\theta \ra 0} \limsup_{N\ra \oo} \sup_\tau \sup_l \mathbb{E}_{\mu^{\scriptscriptstyle{N}}}^N \left[ | \<\tilde{\pi}_{\tau + \theta}^N,G_l\rangle - \<\tilde{\pi}_{\tau}^N,G_l\rangle | \right] =0.
$$
Therefore, we just have to prove that
$$
\limsup_{N\ra \oo} \sup_\tau \mathbb{E}_{\mu^{\scriptscriptstyle{N}}}^N \left[ \frac{1}{N} \sum_{x=-CN}^{0} \eta_{\tau}^N(x) \right] 
$$
converges to $0$ as $C \ra 0$. The previous expression is dominated by
$$
\limsup_{N\ra \oo} \mathbb{E}_{\mu^{\scriptscriptstyle{N}}}^N \left[ \sup_{0\le s\le T} \frac{1}{N} \sum_{x=-CN}^{0} \eta_s^N(x) \right] \, ,
$$
which, by the coupling also described in the end of the proof of (\ref{eq:decay}), is bounded by 
$$
\limsup_{N\ra \oo} \bar{\mathbb{E}}_{\tilde{\nu}_\alpha}^N \left[ \sup_{0\le s\le T} \frac{1}{N} \sum_{x=-CN}^{0} \eta_s^N(x) \right] \, ,
$$
where $\bar{\mathbb{E}}_{\tilde{\nu}_\alpha}^N$ denotes the expectation with respect to the distribution of the nearest-neighbor, symmetric,  space-homogeneous zero-range process on $\mathbb{Z}_-$ with reflection at the origin speeded up by $N^2$ and with initial measure $\tilde{\nu}_\alpha$. Therefore, we conclude the proof of (\ref{eq:stdecay}) with the following result:
\begin{equation}
\label{eq:empmeasdecay}
\lim_{C \ra 0} \, \limsup_{N\ra \oo} \bar{\mathbb{E}}_{\tilde{\nu}_\alpha}^N \left[ \sup_{0\le s\le T} \frac{1}{N} \sum_{x=-CN}^{0} \eta_s^N(x) \right] = 0 \, .
\end{equation}
To prove this, we fix a smooth positive function $H_C :\mathbb{R}_+ \ra \mathbb{R}$ such that $H_C \equiv 1$ on $[0,C]$ and its support is in $[0,2C]$.
Then the expectation in the previous expression is bounded by
$$
\bar{\mathbb{E}}_{\tilde{\nu}_\alpha}^N \left[ \sup_{0\le s\le T} \<\bar{\pi}_s^N,H_C \rangle \right] \, ,
$$
where $\bar{\pi}$ is the empirical measure associated to the zero-range with reflection. Thus an upper bound for (\ref{eq:empmeasdecay}) is given by
\begin{eqnarray}
\label{eq:partition}
\lefteqn{ \limsup_{K \ra \oo} \, \limsup_{N \ra \oo} \bar{\mathbb{E}}_{\tilde{\nu}_\alpha}^N \left[ \max_{0\le i\le K} \<\bar{\pi}_{\frac{iT}{K}}^N,H_C \rangle \right] + } \nn \\
& & + \limsup_{K \ra \oo} \, \limsup_{N \ra \oo} \bar{\mathbb{E}}_{\tilde{\nu}_\alpha}^N \left[ \sup_{|s-t| \le \frac{T}{K}} | \<\bar{\pi}_t^N,H_C \rangle  - \<\bar{\pi}_s^N,H_C \rangle | \right] \, . 
\end{eqnarray} 
Now, choosing $\beta >0$ sufficiently small such that 
$
E_{\tilde{\nu}_\alpha} \left[ \exp \{ \beta \eta(0) \} \right]  
$ 
is finite, we have that the expectation in the first term of (\ref{eq:partition}) is dominated by
$$
\frac{1}{\beta N} \log \bar{\mathbb{E}}_{\tilde{\nu}_\alpha}^N \left[ \exp \left\{ \beta \max_{0\le i\le K} \sum_{x=-2CN}^{0} \eta_{\frac{iT}{K}}^N(x) \right\} \right] 
$$
which is bounded by
$$
\frac{\log K}{\beta N} + 2 \beta^{-1} C E_{\tilde{\nu}_\alpha} \left[ \exp \{ \beta \eta(0) \} \right],
$$
since $\exp\{\max_{1\le i\le K} a_i\} \le \sum_{1\le i\le K} \exp a_i$ and $\tilde{\nu}_\alpha$ is a product measure invariant for the zero-range reflected at the origin.
Therefore the first term of (\ref{eq:partition}) is of order $O(C)$.
On the other hand, the expectation in the second term of (\ref{eq:partition}) is proved to be of order $O(C^{-1} K^{-1})$ as $N \ra \oo$ by standard techniques, using the martingale associated to the zero-range process in the semi-infinite space with reflexion at the origin (see chapter 5 in \cite{kipnislandim}), which means that 
the second term of (\ref{eq:partition}) is zero. This proves (\ref{eq:empmeasdecay}). $\square$

\medskip
\subsection{The proof of Lemma (\ref{lemma:weaksol}) }
\label{sec:weaksol}

Let $G:[0,T] \times \mathbb{R}_+ \ra \mathbb{R}$ be a smooth function with compact support. By Doob inequality, for every $\delta > 0$,
$$
\mathbb{P}_{\mu^{\scriptscriptstyle{N}}}^N \left[ \sup_{0\le t \le T} |M_t^{G,N}| \ge \delta \right] \le 4 \delta^{-2} \mathbb{E}_{\mu^{\scriptscriptstyle{N}}}^N \left[ (M_T^{G,N})^2 \right]
= 4 \delta^{-2} \mathbb{E}_{\mu^{\scriptscriptstyle{N}}}^N \left[  \<M^{G,N}\rangle_T \right], 
$$
which, by the explicit formula for the quadratic variation of $M_t^{G,N}$, is bounded by
$$
\frac{C(G)}{N} \left( \theta + \mathbb{E}_{\mu^{\scriptscriptstyle{N}}}^N \left[ \int_0^T N g(\eta_s(0)) ds \right] \right).
$$
Thus, by (\ref{eq:decay}), for every $\delta >0$, 
\begin{equation}
\label{eq:supmart}
\lim_{N \ra \oo} \mathbb{P}_{\mu^{\scriptscriptstyle{N}}}^N \left[ \sup_{0\le t \le T} |M_t^{G,N}| \ge \delta \right] =0.
\end{equation} 
Using (\ref{eq:integrand}) to expand the martingale expression in (\ref{eq:martingal}) and since the Taylor expansion gives us that
$$
\sup_{s \in [0,T]}
\left| N [G(s,x+1/N)-G(s,x/N)] - \nabla G(s,x/N) \right| \le \frac{C(G)}{N},
$$
we may replace $\Delta_N$ and $\nabla_N$ in (\ref{eq:supmart})  by the usual laplacian and gradient, i.e.,
\begin{eqnarray}
\lefteqn{\!\!\!\!\!\!\!\!\!\!\!\!\!\!
\lim_{N \ra \oo} \mathbb{Q}_{\mu^N}^N \left[ \sup_{0 \le t \le T} \left| \<\pi_t^N,G\rangle - \<\pi_0^N,G\rangle - \int_0^t  ds \{ \<\pi_s^N,\partial_s G\rangle + \frac{1}{2} \< \pi_s^N,\Delta G \rangle + \right. \right.} \nn \\
& & \qquad \qquad \left. \left. + \frac{1}{2} \nabla G(s,0) \xi_s(1) +  N g(\eta_s(0)) \< \pi_s^N, \nabla G \rangle  \} \right| >\delta \right] =0 \nn
\end{eqnarray}
for all $\delta>0$.

In the previous expression we claim that we can obtain the integral term as a function of the empirical measure replacing $Ng(\eta_s(0))$ by $a_s$, given in (\ref{eq:a}), and $\xi_s(1)$ by a mean of $\xi$ over small boxes around 1. At first we are going to justify the replacement of $Ng(\eta_s(0))$ by $a_s$. This is the content of Lemma \ref{lemma:replacement} just below. Let us note before that for all $\delta>0$, and $0\le t \le T$,
\begin{equation}
\label{eq:replacement}
\lim_{N\ra \oo } \mathbb{P}_{\mu^{\scriptscriptstyle{N}}}^N \left[ \left| \int_0^t ds \{ Ng(\eta_s(0)) - a_s \} \right| > \delta \right] = 0, 
\end{equation}
which indicates that the right candidate to replace $Ng(\eta_s(0))$ is $a_s$. Actually, this follows from (\ref{eq:nghl}), (\ref{eq:martdecay}) and Proposition 5.1 in \cite{landimollavolchan}, which states that for all $\delta>0$, and $0\le t \le T$,
$$
\lim_{N\ra \oo} \tilde{\mathbb{Q}}_{\mu^{\scriptscriptstyle{N}}}^N \left[ \left| \<\tilde{\pi}_t^N,1\rangle - \<\tilde{\pi}_0^N,1\rangle -v_t \right| > \delta \right] =0, 
$$
with, by (\ref{eq:v}) and (\ref{eq:a}),
$$
v_t = \int_0^t a_s ds = \int_0^{\oo} \{\rho(t,u) -\rho(0,u)\} du .
$$
Here, 
$$
\<\tilde{\pi}_t^N,1\rangle - \<\tilde{\pi}_0^N,1\rangle \quad \textrm{and} \quad
\int_0^{\oo} \{\rho(t,u) -\rho(0,u)\} du
$$
are to be understood respectively as 
$$
\lim_{l \ra \oo} \{ \<\tilde{\pi}_t^N,H_l\rangle - \<\tilde{\pi}_0^N,H_l\rangle \} \quad \textrm{and} \quad \lim_{l \ra \oo} \, \int_0^{\oo} H_l(u) \{\rho(t,u) -\rho(0,u)\} du \, ,
$$
where $H_l$ is defined in section \ref{sec:tightness}. 

\medskip
\begin{lemma} 
\label{lemma:replacement} 
For every smooth function $G:\mathbb{R}_+ \ra \mathbb{R}$ with compact support and $\delta >0$
$$
\limsup_{N\ra \oo } \mathbb{Q}_{\mu^{\scriptscriptstyle{N}}}^N \left[ \sup_{0 < t\le T} \left| \int_0^t \{ Ng(\eta_s(0)) - a_s \} \<\pi_s^N, G\rangle ds \right| > \delta \right] = 0 \, .
$$
\end{lemma} 

\no {\textbf{Proof:}} The supremum in the statement is bounded above by
\begin{eqnarray}
\label{replacement}
\lefteqn{ \!\!\!\!\!\!\!\!\!\!\!\!\!\!\!\!\!\!\!\!\!\!\!\!
\max_{0 < j \le K} \left| \int_0^{\frac{jT}{K}} \{ Ng(\eta_s(0)) - a_s \} \<\pi_s^N, G\rangle ds \right| + } \nn \\
& & + \sup_{0 < t\le T} \int_t^{t+\frac{T}{K}} Ng(\eta_s(0)) ds + \frac{C(G)}{K} \sup_{0\le t \le T} a_s \, ,
\end{eqnarray}
for every $K>0$. The third term clearly goes to $0$ as $K$ goes to $\oo$. To deal with the second term we are going to show that     
\begin{equation}
\label{eq:equicont}
\limsup_{\theta \ra 0} \limsup_{N \ra \oo} \mathbb{E}_{\mu^{\scriptscriptstyle{N}}}^N \left[ \sup_{0 < t\le T} \int_t^{t+\theta} Ng(\eta_s(0)) ds \right] = 0 .
\end{equation}
To show (\ref{eq:equicont}), recall formula (\ref{eq:nghl2}) with the random time $\tau$ replaced by $t$.
By Fatou's Lemma, it is enough to prove that
$$
\lim_{\theta \ra 0} \limsup_{N \ra \oo} \limsup_{l \ra \oo} \mathbb{E}_{\mu^{\scriptscriptstyle{N}}}^N \left[ \sup_{0 < t\le T} | \tilde{M}^{H_l,N}_{t + \theta} - \tilde{M}^{H_l,N}_t | \right] = 0 .
$$
and that
$$
\lim_{\theta \ra 0} \limsup_{N \ra \oo} \limsup_{l \ra \oo} \mathbb{E}_{\mu^{\scriptscriptstyle{N}}}^N \left[ \sup_{0 < t\le T} | \<\tilde{\pi}^N_{t + \theta},H_l\rangle - \<\tilde{\pi}^N_{t},H_l\rangle | \right] = 0 .
$$
Considering the supremum out of the expectation in the last two expression, both of them were proved in the last section. The same proof carried out there can be applied for the case where the supremum is inside the expectation, since we have (\ref{eq:empmeasdecay}) and the expectation of the supremum of a martingale is dominated by the expectation of its quadratic variation. 

\smallskip 
So it remains to consider the first term in (\ref{replacement}). Since we first make $N \ra \oo$ and then $K \ra \oo$, we only have to show that
\begin{equation}
\label{replacement2}
\limsup_{N\ra \oo } \mathbb{Q}_{\mu^{\scriptscriptstyle{N}}}^N \left[ \left| \int_0^t \{ Ng(\eta_s(0)) - a_s \} \<\pi_s^N, G\rangle ds \right| > \delta \right] = 0 \, .
\end{equation}
for all $0 < t \le T$. Fix $C>0$ and let $\mathcal{P} : \, 0=t_0 < t_1< \ldots <t_n=t$ be a partition on $[0,t]$ such that $|\mathcal{P}| := \max(t_{i+1}-t_i) \le C \min(t_{i+1}-t_i)$. Bound (\ref{replacement2}) by
\begin{eqnarray}
\label{eq:decomp}
\lefteqn{\sum_{i=0}^{n-1} \left\{ \left| \int_{t_i}^{t_{i+1}} Ng(\eta_s(0)) \{ \<\pi_s^N,  G\rangle - \<\pi_{t_i}^N, G\rangle \} ds \right| + \right. } \nn \\ 
& & \qquad  + \left| \int_{t_i}^{t_{i+1}} a_s \{ \<\pi_s^N, G\rangle - \<\pi_{t_i}^N, G\rangle \} ds \right| + \\ 
& & \qquad + |\<\pi_{t_i}^N, G\rangle | \left. \left| \int_{t_i}^{t_{i+1}} \{Ng(\eta_s(0)) - a_s \} ds \right| \right\} . \nn 
\end{eqnarray}
In order to prove (\ref{replacement2}), we start estimating the $\mathbb{Q}_{\mu^{\scriptscriptstyle{N}}}^N$ probability of the first term in (\ref{eq:decomp}) to be greater than $\delta$. Using (\ref{eq:martingal}) and (\ref{eq:integrand}), we obtain an upper bound of 
\begin{eqnarray}
\label{eq:primeirotermo}
\lefteqn{ \sum_{i=0}^{n-1} \int_{t_i}^{t_{i+1}} Ng(\eta_s(0)) \left\{ | M_s^{G,N} - M_{t_i}^{G,N}| + \right. } \nn \\
& & \qquad \left. + C(G) (t_{i+1} - t_i) + C(G) \int_{t_i}^{t_{i+1}} Ng(\eta_s(0)) ds \right\} ds \, .
\end{eqnarray}
Thus, we are going consider separately each term in the previous expression:

\smallskip
\no \textbf{Claim 1:}  
$$
\lim_{|\mathcal{P}| \ra 0} \lim_{N \ra \oo} \mathbb{E}_{\mu^{\scriptscriptstyle{N}}}^N \left[  \sum_{i=0}^{n-1} \int_{t_i}^{t_{i+1}} Ng(\eta_s(0)) | M_s^{G,N} - M_{t_i}^{G,N}| ds \right] = 0 \, . 
$$
\no \textbf{Proof of Claim 1:} By Cauchy-Schwarz inequality the expectation in the statement is bounded above by
\begin{equation}
\label{eq:claim1a}
N \sum_{i=0}^{n-1} \mathbb{E}_{\mu^{\scriptscriptstyle{N}}}^N \left[ \int_{t_i}^{t_{i+1}} g(\eta_s(0)) ds \right]^{\frac{1}{2}} \,  \mathbb{E}_{\mu^{\scriptscriptstyle{N}}}^N \left[ \int_{t_i}^{t_{i+1}} | M_s^{G,N} - M_{t_i}^{G,N}|^2 ds \right]^{\frac{1}{2}}  
\end{equation}
From (\ref{eq:quavar}) and (\ref{eq:aldous2}) we obtain the following estimate
\begin{eqnarray}
\lefteqn{ \!\!\!\!\!\!\!\!\!\!\!
\mathbb{E}_{\mu^{\scriptscriptstyle{N}}}^N \left[ \int_{t_i}^{t_{i+1}} | M_s^{G,N} - M_{t_i}^{G,N}|^2 ds \right] = \int_{t_i}^{t_{i+1}} \mathbb{E}_{\mu^{\scriptscriptstyle{N}}}^N \left[ \<M^{G,N}\rangle_s - \<M^{G,N}\rangle_{t_i} \right] ds } \nn \\ 
& & \quad \le C(G) \left\{ \frac{(t_{i+1}-t_i)^2}{N} + (t_{i+1}-t_i) \mathbb{E}_{\mu^{\scriptscriptstyle{N}}}^N \left[ \int_{t_i}^{t_{i+1}} g(\eta_s(0)) ds \right] \right\} .  \nn
\end{eqnarray}
Therefore, since $(a+b)^{\frac{1}{2}} \le a^{\frac{1}{2}} + b^{\frac{1}{2}}$, (\ref{eq:claim1a}) is bounded by
\begin{eqnarray}
\label{eq:claim1b}
\lefteqn{C(G) \sum_{i=0}^{n-1} (t_{j+1} - t_j) \mathbb{E}_{\mu^{\scriptscriptstyle{N}}}^N \left[ \int_{t_i}^{t_{i+1}} N g(\eta_s(0)) ds \right]^{\frac{1}{2}} + } \nn \\
& & \qquad \qquad + C(G) \sum_{i=0}^{n-1} (t_{j+1} - t_j)^{\frac{1}{2}} \mathbb{E}_{\mu^{\scriptscriptstyle{N}}}^N \left[ \int_{t_i}^{t_{i+1}} N g(\eta_s(0)) ds \right] .
\end{eqnarray}
Applying Schwarz inequality to the first term, we have that (\ref{eq:claim1b}) is dominated by
$$
C(G) \max(t_{j+1} - t_j)^{\frac{1}{2}} \left\{ 1 + \mathbb{E}_{\mu^{\scriptscriptstyle{N}}}^N \left[ \int_0^t N g(\eta_s(0)) ds \right]^{\frac{1}{2}} \right\}^2 .
$$
Together with (\ref{eq:decay}) this proves Claim 1. $\square$

\smallskip

\no \textbf{Claim 2:}
$$
\lim_{|\mathcal{P}|\ra 0} \lim_{N \ra \oo} \mathbb{E}_{\mu^{\scriptscriptstyle{N}}}^N \left[  \sum_{i=0}^{n-1} (t_{i+1} - t_i) \int_{t_i}^{t_{i+1}} Ng(\eta_s(0)) ds \right] = 0 \, . 
$$
\no \textbf{Proof of Claim 2:} This expectation is bounded by
$$
\max(t_{j+1}-t_j) \mathbb{E}_{\mu^{\scriptscriptstyle{N}}}^N \left[ \int_0^t N g(\eta_s(0)) ds \right] .
$$
Together with (\ref{eq:decay}) this proves Claim 2. $\square$

\smallskip

\no \textbf{Claim 3:} For every $\delta>0$
$$
\lim_{|\mathcal{P}|\ra 0} \lim_{N \ra \oo} \mathbb{P}_{\mu^{\scriptscriptstyle{N}}}^N \left[  \sum_{i=0}^{n-1} \left\{ \int_{t_i}^{t_{i+1}} Ng(\eta_s(0)) ds \right\}^2 > \delta \right] = 0 \, . 
$$
\no \textbf{Proof of Claim 3:} The sum in the previous expression is bounded by
$$
2 \sum_{i=0}^{n-1} \left| \int_{t_i}^{t_{i+1}} \left[ Ng(\eta_s(0)) - a_s \right] ds \right|^2 + 2 \sum_{i=0}^{n-1}  \left\{ \int_{t_i}^{t_{i+1}} a_s ds \right\}^2 .
$$
By (\ref{eq:replacement}) the first term in this last expression goes to $0$ in probability as $N \ra \oo$, while the second term converges to $0$ as $|\mathcal{P}|\ra 0$. Hence, Claim 3 holds. $\square$

\smallskip
Since (\ref{eq:primeirotermo}) is an upper bound for the first term in (\ref{eq:decomp}), from claim 1-3 we conclude that for all $\delta>0$ 
$$
\lim_{N \ra \oo} \mathbb{Q}_{\mu^{\scriptscriptstyle{N}}}^N \left[ \sum_{i=0}^{n-1} \left| \int_{t_i}^{t_{i+1}} Ng(\eta_s(0)) \{ \<\pi_s^N, G\rangle - \<\pi_{t_i}^N, G\rangle \} ds \right| > \delta \right] =0 .
$$

Analogously, considering the second term in (\ref{eq:decomp}), we show that for all $\delta >0$
$$
\lim_{N \ra \oo} \mathbb{Q}_{\mu^{\scriptscriptstyle{N}}}^N \left[ \sum_{i=0}^{n-1} \left| \int_{t_i}^{t_{i+1}} a_s \{ \<\pi_s^N, G\rangle - \<\pi_{t_i}^N, G\rangle \} ds \right| > \delta \right] = 0 .
$$

It remains to consider the last term in (\ref{eq:decomp}), but again it is bounded by
$$
C(G) \sum_{i=0}^{n-1} \left| \int_{t_i}^{t_{i+1}} \left[ Ng(\eta_s(0)) - a_s \right] ds \right|
$$
and by (\ref{eq:replacement}) it converges to $0$ in probability as $N \ra \oo$.   This concludes the proof. $\square$ 

\bigskip

Now we consider the replacement of $\xi_s(1)$. At this point we need the condition (E5) described in section \ref{sec:initialmeasures}. In section \ref{sec:rementcond} we justify how this condition can be removed from the hypotheses to prove Theorem \ref{theorem:hbep}.

\bigskip
\begin{lemma}
\label{lemma:replacementbymean}
For a sequence of initial measures $\{\mu^N, N>1\}$ satisfying (E1)-(E5), we have for every continuously differentiable function $H:[0,T] \ra \mathbb{R}$ that
$$
\limsup_{\epsilon \ra 0} \limsup_{N \ra \oo} \mathbb{E}_{\mu^{\scriptscriptstyle{N}}}^N \left[ \sup_{0\le t \le T} \left| \int_0^t H(s) \Big\{ \xi_s(1) - \frac{1}{\epsilon N} \sum_{x=1}^{[\epsilon N]} \xi_s(x) \Big\} ds \right| \right] = 0 \, .
$$
\end{lemma}

\no {\textbf{Proof:}} As in the proof of Lemma \ref{lemma:replacement}, since 
$$
\lim_{\theta \ra 0} \limsup_{\epsilon \ra 0} \limsup_{N \ra \oo} \mathbb{E}_{\mu^{\scriptscriptstyle{N}}}^N \left[ \sup_{0\le t \le T} \int_t^{t+\theta} H(s) \Big\{ \xi_s(1) - \frac{1}{\epsilon N} \sum_{x=1}^{[\epsilon N]} \xi_s(x) \Big\} ds \right]=0,
$$
we may omit the supremum in the statement. It is also easily seen that we can replace $\xi_s(0)$ by $\frac{[\epsilon N]}{\epsilon N} \xi_s(0)$ in the statement and it is enough to show that
\begin{equation}
\label{eq:expU}
\limsup_{\epsilon \ra 0} \limsup_{N \ra \oo} \mathbb{E}_{\mu^{\scriptscriptstyle{N}}}^N \left[ \left| \int_0^t U_{\epsilon}^{N}(s,\xi_s) ds \right| \right] ,
\end{equation}
where $U_{\epsilon}^{N}(s,\xi)=H(s)V_{\epsilon}^{N}(\xi)$, $0\le s \le t$, with
$$
V_{\epsilon}^{N}(\xi) = \frac{1}{\epsilon N} \sum_{x=1}^{[\epsilon N]} ( \xi(0) - \xi(x) ) \, .
$$
By the entropy inequality we have that the expectation in (\ref{eq:expU}) is bounded above by
\begin{equation}
\label{eq:entest}
\frac{ H(\mu^{\scriptscriptstyle{N}}|\nu_{\alpha,\gamma}^N) }{A_N} + \frac{1}{A_N} \log \mathbb{E}_{\nu_{\alpha,\gamma}^N}^N \left[ \exp \left\{ \left| \int_0^t A_N U_{\epsilon}^{N}(s,\xi_s) ds \right| \right\} \right] 
\end{equation}
for every $A_N$, where $\nu_{\alpha,\gamma}^N$ is defined in section \ref{sec:initialmeasures}.. 

We are going to estimate the second term in (\ref{eq:entest}). Since $e^{|x|} \le e^{x} + e^{-x}$ and $\limsup_{N} N^{-1} \log \{ a_N + b_N \} \le \max \{ \limsup_{N} N^{-1} \log a_N , \limsup_{N} N^{-1} \log b_N \}$, we may suppress the absolute value in the exponent. Define
$$
(P_{s,t}^N f)(\eta,\xi) = \mathbb{E}_{(\eta,\xi)}^N \left[ 
f(\eta_{t-s} , \xi_{t-s}) \exp \left\{ \int_0^{t-s} A_N U(s+r, \xi_r) dr \right\}  \right]
$$
for every bounded function $f$ on $\Omega$. We have that
$$
\mathbb{E}_{\nu_{\alpha,\gamma}^N}^N \left[ \exp \left\{ \int_0^t A_N U_{\epsilon}^{N}(s,\xi_s) ds \right\} \right] = \int P_{0,t}^N 1 \ d\nu_{\alpha,\gamma}^N \le \left\{ \int (P_{0,t}1^N)^2 \ d\nu_{\alpha,\gamma}^N \right\}^{\frac{1}{2}} \, .
$$
In order to obtain an upper bound for the right hand term in this inequality we will show below, finishing this section, that 
\begin{eqnarray}
\label{eq:derPst}
\!\!\!\!\!\!\!\!\!\!\!\!\!\!\!\!\!\!\!
- \frac{1}{2} \partial_s  \int (P_{s,t}^N 1)^2  d\nu_{\alpha,\gamma}^N  & \le &
\left\{ \frac{\beta N}{1-\alpha} + \left( \frac{B A_N}{2} - N^2 \right) \sup_f D(f) \right. +  \nn \\
& & \quad \qquad \left. + \frac{\epsilon N A_N}{B} \| H(u)\|_{\oo}^2 \right\}\int (P_{s,t}^N 1)^2  d\nu_{\alpha,\gamma}^N \, , 
\end{eqnarray}
for every $B>0$, where the supremum in $f$ is taken over all densities with respect to $\nu_\alpha$, $\| \cdot \|_{\oo}$ denote the supremum norm and $D(f)$ denote the Dirichlet form
$$
- \int \sqrt{f} L \sqrt{f} d\nu_\alpha .
$$
As usual, we can divide both sides of (\ref{eq:derPst}) by $\int (P_{s,t}^N 1)^2  d\nu_{\alpha,\gamma}^N > 0$ and integrate over $[0,t]$ to obtain that
$$
\int (P_{s,t}^N 1)^2  d\nu_{\alpha,\gamma}^N \le \exp \Big\{ 2t \Big( \frac{\beta N}{1-\alpha} + \Big( \frac{B A_N}{2} - N^2 \Big) \sup_f D(f) + \frac{\epsilon N A_N}{B} \| H(u)\|_{\oo}^2 \Big) \Big\}
$$
In particular, since by (\ref{eq:entbound}) and (E5) we have that $H(\mu^N | \nu_{\alpha,\gamma}^N) \le CN$, for some $C>0$, it follows from the previous inequality that (\ref{eq:entest}) is bounded by
$$
\frac{CN}{A_N} + \left( \frac{\beta N}{(1-\alpha) A_N}  + \left( \frac{B}{2} - \frac{N^2}{A_N} \right) \sup_f D(f) + \frac{\epsilon N}{B} \| H(u)\|_{\oo}^2 \right) t \, .
$$

Now, choosing $B=2 \sqrt{\epsilon}N$ and $A_N=N/\sqrt{\epsilon}$, it turns out that 
$$
\mathbb{E}_{\mu^{\scriptscriptstyle{N}}}^N \left[ \left| \int_0^t U_{\epsilon}^{N}(s,\xi_s) ds \right| \right] \le \{C + \beta (1-\alpha)^{-1} + 2^{-1} \| H(u)\|_{\oo}^2 \} \sqrt{\epsilon} T \, ,
$$
which goes to $0$ as $\epsilon \ra \oo$, proving the Lemma. $\square$

\medskip

\no \textbf{Proof of (\ref{eq:derPst}):} By Feynman-Kac formula, $\{ P_{s,t}: 0\le s\le t\}$ is a semigroup of operators associated to the non-homogeneous generator $\mathcal{L}_s = N^2 \mathcal{L} + A_N U_{\epsilon}^{N} (s,\cdot)$. Moreover, the first Chapman-Kolmogorov equation holds: $\partial_s P_{s,t} = - \mathcal{L}_s P_{s,t}$. Hence
\begin{eqnarray}
\lefteqn{ - \frac{1}{2} \partial_s  \int (P_{s,t}^N 1)^2 \, d\nu_{\alpha,\gamma}^N = \int \mathcal{L}_s (P_{s,t}^N 1) \, P_{s,t}^N 1 \  d\nu_{\alpha,\gamma}^N } \nn \\
& & = \int N^2 \tilde{L} (P_{s,t}^N 1) \, P_{s,t}^N 1 \ d\nu_{\alpha,\gamma}^N + \int N^2 L_b (P_{s,t}^N 1) \, P_{s,t}^N 1 \ d\nu_{\alpha,\gamma}^N + \nn \\
& & \qquad + \int \{ N^2 L + A_N U_{\epsilon}^{N} (s,\cdot)\} (P_{s,t}^N 1) \, P_{s,t}^N 1 \ d\nu_{\alpha,\gamma}^N \, . \nn
\end{eqnarray}
We shall estimate separately each term in this expression.

\medskip
\no \textbf{Claim 1:} 
$$
\int \tilde{L} (P_{s,t}^N 1) \, P_{s,t}^N 1 \ d\nu_{\alpha,\gamma}^N \le \frac{1}{2} \int  g(\eta(0)) \, (P_{s,t}^N 1)^2 \, d\nu_{\alpha,\gamma}^N \, .
$$
\no \textbf{Proof of Claim 1:} Denote $P_{s,t}^N 1$ by h. We have that
$$
\int (\tilde{L} h) \, h \ d\nu_{\alpha,\gamma}^N = \sum_{x\le -1} \int \{ \tilde{L}_{x,x+1}h + \tilde{L}_{x+1,x}h \} \, h \ d\nu_{\alpha,\gamma}^N \, .
$$
Recall from section \ref{sec:initialmeasures} the definition of $\nu^N_{\alpha,\gamma}$ and $\gamma_x^N$, $x\le 0$, $N\ge 1$. Considering the change of variables $\sigma^{x,x+1} \eta = \tilde{\eta}$ we have that
$$
\frac{d\nu_{\alpha,\gamma}^N(\eta)}{d\nu_{\alpha,\gamma}^N(\tilde{\eta})}
= \frac{\gamma_{x}^N \, g(\tilde{\eta}(x+1))}{\gamma_{x+1}^N \, g(\eta(x))}
$$
and therefore
$$
 \int \{ \tilde{L}_{x,x+1}h + \tilde{L}_{x+1,x}h \} \, h \ d\nu_{\alpha,\gamma}^N
$$
can be written as
\begin{eqnarray}
\lefteqn{ - \int \frac{g(\eta(x))}{2} [h(\sigma^{x,x+1}\eta, \xi) - h(\eta,\xi)]^2 \, d\nu_{\alpha,\gamma}^N } \nn \\
& & - \int \frac{g(\eta(x+1))}{2} [h(\sigma^{x+1,x} \eta, \xi) - h(\eta,\xi)]^2 \, d\nu_{\alpha,\gamma}^N  \nn \\
& & + \frac{1}{2} \left( \frac{\gamma_{x+1}^N}{\gamma_x^N} - 1 \right)\int g(\eta(x)) \, h(\eta,\xi)^2 \, d\nu_{\alpha,\gamma}^N \nn \\
& & + \frac{1}{2} \left( \frac{\gamma_x^N}{\gamma_{x+1}^N} - 1 \right)\int g(\eta(x+1)) \, h(\eta,\xi)^2 \, d\nu_{\alpha,\gamma}^N \, . \nn 
\end{eqnarray}
In the previous summation we can neglect the first two terms, which are negative, and add the last two terms in $x$, for $x \le -1$, obtaining that 
\begin{eqnarray}
\lefteqn{ \!\!\!\!\!\!\!\!\!\!\!\!\!\!\!\!\!\!
\int (\tilde{L} h) \, h \ d\nu_{\alpha,\gamma}^N \le \frac{1}{2}  \left( \frac{\gamma_{-1}^N}{\gamma_0^N} - 1 \right)\int g(\eta(0)) \, h(\eta,\xi)^2 \, d\nu_{\alpha,\gamma}^N + } \nn \\
& & \qquad + \frac{1}{2} \sum_{x\le -1} \frac{\Delta \gamma^N (x)}{\gamma_x^N} \int g(\eta(x)) \, h(\eta,\xi)^2 \, d\nu_{\alpha,\gamma}^N \, ,\nn
\end{eqnarray}
where $\Delta \gamma^N (x) = \gamma_{x+1}^N + \gamma_{x-1}^N - 2 \gamma_x^N$. Observing that $\Delta \gamma^N (x)$ is zero except at $x=-N+1$, when it is negative, and that $\gamma_{-1}^N / \gamma_0^N = 2$, we have shown Claim 1. $\square$

\medskip
\no \textbf{Claim 2:} 
$$
\int L_b (P_{s,t}^N 1) \, P_{s,t}^N 1 \ d\nu_{\alpha,\gamma}^N \le \frac{\beta}{2 (1-\alpha) N} \int (P_{s,t}^N 1)^2 \, d\nu_{\alpha,\gamma}^N
- \frac{1}{2} \int  g(\eta(0)) \, (P_{s,t}^N 1)^2 \, d\nu_{\alpha,\gamma}^N \, .
$$
\no \textbf{Proof of Claim 2:} Denote $P_{s,t}^N 1$ by h. We are considering an integral of the form
$$
\int g(\eta(0)) [h(\eta-\varrho_0, \tau \xi) - h(\eta,\xi)] \, h(\eta,\xi) \ d\nu_{\alpha,\gamma}^N \, ,
$$
where $\varrho_0$ and $\tau$ are defined in section \ref{sec:system}. Add and subtract to this the term
$$
\frac{1}{2} \int g(\eta(0)) h^2(\eta-\varrho_0, \tau \xi) \ d\nu_{\alpha,\gamma}^N = \frac{\gamma_0^N}{2 (1-\alpha)} \int (1-\xi(0)) h^2 \  d\nu_{\alpha,\gamma}^N \, .
$$
We end up with three terms, the first is 
$$
- \frac{1}{2} \int g(\eta(0)) [h(\eta-\varrho_0, \tau \xi) - h(\eta,\xi)]^2 \ d\nu_{\alpha,\gamma}^N \, ,
$$
which is negative and may be neglected, and the others are
$$
\frac{\gamma_0^N}{2 (1-\alpha)} \int (1-\xi(0)) h^2 \ d\nu_{\alpha,\gamma}^N - \frac{1}{2} \int  g(\eta(0)) \, h^2 \, d\nu_{\alpha,\gamma}^N
$$
Since $\gamma_0^N = \beta/N$, we have Claim 2. $\square$

\medskip
\no \textbf{Claim 3:} 
\begin{eqnarray}
\lefteqn{\int \{ N^2L + A_N U_{\epsilon}^{N} (s,\cdot)\} (P_{s,t}^N 1) \, P_{s,t}^N 1 \ d\nu_{\alpha,\gamma}^N \le } \nn \\
& & \le \left\{ \sup_f \left\{ \left( \frac{B A_N}{2} - N^2 \right) D(f) \right\} + \epsilon \frac{N A_N}{B} \| H(u)\|_{\oo}^2 \right\} \int (P_{s,t}^N 1)^2 \, d\nu_{\alpha,\gamma}^N \, . \nn
\end{eqnarray}
\no \textbf{Proof of Claim 3:} We have the bound 
$$
\int \{ N^2L + A_N U_{\epsilon}^{N} (s,\cdot)\} (P_{s,t}^N 1) \, P_{s,t}^N 1 \ d\nu_{\alpha,\gamma}^N
\le \Gamma_s^N \int (P_{s,t}^N 1)^2 \ d\nu_{\alpha,\gamma}^N \, ,
$$
where $\Gamma_s^N$ is the greatest eigenvalue of the generator $N^2 L + A_N U_{\epsilon}^{N}(s, \cdot )$. By the variational formula (see Appendix 3 in \cite{kipnislandim}) $\Gamma_s^N$ is equal to
\begin{equation}
\label{eq:varform}
\sup_f \left\{ \int A_N U_{\epsilon}^{N}(s,\xi) f(\xi) \nu_{\alpha}(d\xi) - N^2 D(f) \right\}.
\end{equation}
If we can show that
\begin{equation}
\label{eq:intvf}
\int V^N_\epsilon (\xi)f(\xi) \nu_\alpha(d\xi) \le B D(f) + \frac{2 \epsilon N}{B},
\end{equation}
for any $B>0$, then replacing $B$ by $B/H(s)$ we conclude the proof of the Claim from (\ref{eq:varform}). The left hand side in (\ref{eq:intvf}) is equal to
$$
\frac{1}{\epsilon N} \sum_{x=1}^{[\epsilon N]} \int [\xi(1) - \xi(x)] f(\xi) \nu_\alpha(d\xi)  = \frac{1}{\epsilon N} \sum_{x=1}^{[\epsilon N]} \sum_{y=1}^{x-1} \int [\xi(y) - \xi(y+1)] f(\xi), \nu_\alpha(d\xi)
$$
which can be rewritten as
\begin{eqnarray}
\label{eq:expansao}
& & \frac{1}{\epsilon N} \sum_{x=1}^{[\epsilon N]} \sum_{y=1}^{x-1} \left\{ \int \xi(y)[1 - \xi(y+1)] f(\xi) \nu_\alpha(d\xi) - \right. \nn \\ 
& & \qquad\qquad\qquad\qquad - \left. \int \xi(y+1)[1 - \xi(y)] f(\xi) \nu_\alpha(d\xi) \right\} \nn \\
& & \quad = \frac{1}{\epsilon N} \sum_{x=1}^{[\epsilon N]} \sum_{y=1}^{x-1} \int \xi(y+1)[1 - \xi(y)] \left\{ f(\xi^{y,y+1}) - f(\xi) \right\} \nu_\alpha(d\xi) \, .
\end{eqnarray}
Thus, writing $f(\xi^{y,y+1}) - f(\xi)$ as $\{\sqrt{f(\xi^{y,y+1})} - \sqrt{f(\xi)}\} \, \{\sqrt{f(\xi^{y,y+1})} + \sqrt{f(\xi)}\}$ and applying the elementary inequality $2ab \le B a^2 + B^{-1} b^2$, that holds for every $a,b$ in $\mathbb{R}$ and $B>0$, we have that previous expression is bounded by
\begin{eqnarray}
\lefteqn{ \frac{1}{\epsilon N} \sum_{x=1}^{[\epsilon N]} \sum_{y=1}^{x-1} \left\{ \frac{B}{2} \int \xi(y+1)[1 - \xi(y)] \left\{ \sqrt{f(\xi^{y,y+1})} - \sqrt{f(\xi)} \right\}^2 \nu_\alpha(d\xi) + \right. } \nn \\
& & \quad \qquad \qquad + \left. \frac{B^{-1}}{2} \int \xi(y+1)[1 - \xi(y)] \left\{ \sqrt{f(\xi^{y,y+1})} + \sqrt{f(\xi)} \right\}^2 \nu_\alpha(d\xi) \right\} \nn \\
& &  \le \frac{B}{2} D(f) + \frac{B^{-1}}{\epsilon N} \sum_{x=1}^{[\epsilon N]} \sum_{y=1}^{x-1}  \int \left\{ f(\xi^{y,y+1}) + f(\xi) \right\} \nu_\alpha(d\xi) \le \frac{B}{2} D(f) + \frac{\epsilon N}{B} \, . \nn
\end{eqnarray}
This show (\ref{eq:intvf}). $\square$
 
\smallskip 
Now, is easy to see that (\ref{eq:derPst}) is a consequence of Claim 1-3. Therefore we have proved the Lemma.  $\square$

\medskip
\subsection{An energy estimate}
\label{sec:energyest}

The next result justifies an integration by parts in the expression inside the probability in the statement of Lemma \ref{lemma:weaksol}, proving Theorem \ref{theorem:hbep} under condition (E5).

\medskip
\begin{theorem}
\label{theorem:orderofint} 
Every limit point of the sequence $\mathbb{Q}^N$ is concentrated on paths $\zeta(t,u)du$ with the property that $\zeta(t,u)$ is absolutely continuous whose derivative $\partial_u \zeta(t,u)$ is in $L^2([0,T] \times \mathbb{R}_+)$. Moreover
\begin{eqnarray}
\lefteqn{\int_0^T ds \int_{\mathbb{R}_+} du H(s,u) \partial_u \zeta(s,u) = } \nn \\
& = & - \, \int_0^T ds \left\{ \int_{\mathbb{R}_+} du \partial_u H(s,u) \zeta(s,u) + H(s,0) \lim_{\epsilon \ra 0} \epsilon^{-1} \int_0^{\epsilon} \zeta(s,u)du \right\} \nn
\end{eqnarray} 
for all smooth functions $H:[0,T] \times \mathbb{R}_+ \ra \mathbb{R}$ with compact support.
\end{theorem}

\bigskip
Denote by $C^{0,1}_K ([0,T] \times \mathbb{R}_+)$ the space of continuous functions with compact support on $[0,T] \times \mathbb{R}_+$ which are continuously differentiable in the second variable and consider this space endowed with the norm 
$$
\|H\|_{0,1} = \sum_{n=0}^{\oo} 2^n \{ \|H \, \mathbf{1}\{(n,n+1)\}\|_\oo + \|\partial_u H \, \mathbf{1}\{(n,n+1)\}\|_\oo \} \,.
$$ 
To prove the previous theorem we make use of the following energy estimate:

\medskip
\begin{lemma}
\label{lemma:energyestimate}
There exists $K>0$ such that if $\mathbb{Q}^*$ is a limit point of the sequence $\mathbb{Q}^N$ then
\begin{eqnarray}
\mathbb{E}_{\mathbb{Q}^*} \left[ \sup_H \left\{ \int_0^T ds \left\{ \int_{\mathbb{R}_+} \right. \right. \right. \! \! \! \! \! \! \! \!  & \! &  \! \! \! \! \! \left.  du \partial_u H(s,u) \zeta(s,u) + H(s,0) \lim_{\epsilon \ra 0} \epsilon^{-1} \int_0^{\epsilon} \zeta(s,u)du \right\} -  \nn \\
& \! & - \left. \left. 2 \int_0^T ds \int_{\mathbb{R}_+} du H(s,u)^2 \zeta(s,u) \right\} \right] \le K \, , \nn
\end{eqnarray}
where the supremum is taken over all functions $H$ in $C^{0,1}_K ([0,T] \times \mathbb{R}_+)$.
\end{lemma}
\no \textbf{Proof:} For every $\epsilon>0$, $\delta>0$, $H:\mathbb{R}_+ \ra \mathbb{R}$ smooth function with compact support and $\xi \in \{0,1\}^{\mathbb{Z}^*_+}$, denote by $W_N(\epsilon,\delta,H,\xi)$ the following expression
$$
\sum_{x=1}^\oo H(x/N) \frac{1}{\epsilon N} \{\xi^{\delta N}(x) - \xi^{\delta N}(x+[\epsilon N])\} - \frac{2}{N} \sum_{x=1}^\oo H(x/N)^2 \frac{1}{\epsilon N} \sum_{y=0}^{[\epsilon N]} \xi^{\delta N}(x+y) \, ,
$$
where
$
\xi^{\delta}(x) = \delta^{-1} \sum_{y=x}^{x+\delta} \xi(y) .
$
We claim that there exists $K>0$ such that for any dense subset $\{H_l:l\ge 1\}$ of $C^{0,1}_K ([0,T] \times \mathbb{R}_+)$,
\begin{equation}
\label{eq:A}
\lim_{\delta \ra 0} \lim_{N \ra \oo} \mathbb{E}_{\mu^{\scriptscriptstyle{N}}}^N \left[ \max_{1\le i\le k} \left\{ \int_0^T ds W_N(\epsilon,\delta,H_i(s,\cdot),\xi_s) \right\} \right] \le K \, ,
\end{equation}
for every $k\ge 1$ and every $\epsilon >0$. We postpone the proof of (\ref{eq:A}), using it, since $\mathbb{Q}^*$ is a weak limit point of the sequence $\mathbb{Q}_N$, it follows that
\begin{eqnarray}
\label{eq:A1}
\lefteqn{\limsup_{\delta \ra 0} E_{\mathbb{Q}^*} \left[ \max_{1\le i\le k} \left\{  \int_0^T ds \int_{\mathbb{R}_+} du \right. \right. } \nn \\
& & \left\{ \epsilon^{-1} H_i(s,u) \left( \delta^{-1} \int_u^{u+\delta} \zeta_s \, dv - \delta^{-1} \int_{u + \epsilon}^{u+\epsilon+\delta} \zeta_s \, dv \right) \right. - \\ 
& & \qquad \quad \left. \left. \left. - 2 \epsilon^{-1} H_i(s,u)^2 \int_u^{u+\epsilon} dv \left( \delta^{-1} \int_v^{v+\delta} \zeta_s \, dv^\prime \right) \right\} \right\} \right] \le K \, , \nn
\end{eqnarray}
for every $k\ge 1$. Since,
$$
\epsilon^{-1} \int_{\mathbb{R}_+} du H(u) \left\{ \delta^{-1} \int_u^{u+\delta} \zeta_s \, dv - \delta^{-1} \int_{u+\epsilon}^{u+\epsilon+\delta} \zeta_s \, dv \right\} 
$$
is equal to
$$
\int_\epsilon^\oo du \left\{ \frac{H(u)-H(u-\epsilon)}{\epsilon} \right\} \left\{ \delta^{-1} \int_u^{u+\delta} \zeta_s \, dv \right\} + \epsilon^{-1} \int_0^\epsilon du H(u) \left\{ \delta^{-1} \int_u^{u+\delta} \zeta_s \, dv \right\} \, ,
$$
letting $\delta \ra 0$ and then $\epsilon \ra 0$, it follows from (\ref{eq:A1}) that
\begin{eqnarray}
E_{\mathbb{Q}^*} \left[ \max_{1\le i \le k} \left\{ \int_0^T ds \int_{\mathbb{R}_+} du \left\{ \partial_uH_i(s,u) \zeta(s,u) -2 H_i(s,u)^2 \zeta(s,u) \right\} \right. \right. + \qquad & & \nn \\
+ \left. \left. H(s,0) \lim_{\epsilon \ra 0} \epsilon^{-1} \int_0^\epsilon \zeta(s,u) du \right\} \right] \le K \, . & &\nn
\end{eqnarray}
To conclude the proof we apply the monotone convergence theorem, noting that
\begin{eqnarray}
\int_0^T ds \int_{\mathbb{R}_+} du \left\{ \partial_uH_i(s,u) \zeta(s,u) -2 H_i(s,u)^2 \zeta(s,u) \right\}  + \qquad & & \nn \\
+ H(s,0) \lim_{\epsilon \ra 0} \epsilon^{-1} \int_0^\epsilon \zeta(s,u) du \, ,  & & \nn
\end{eqnarray}
is continuous as a real function on $C^{0,1}_K ([0,T] \times \mathbb{R}_+)$. $\square$

\smallskip

\no \textbf{Proof of (\ref{eq:A})}. Since $H$ is a continuous function, an integration by parts justify the replacement of $W_N(\epsilon,\delta,H,\xi)$ as $\delta \ra \oo$ in (\ref{eq:A}) by
\begin{equation}
\label{eq:W}
\sum_{x=1}^\oo H(x/N) \frac{1}{\epsilon N} \{\xi(x) - \xi(x+[\epsilon N])\} - \frac{2}{N} \sum_{x=1}^\oo H(x/N)^2 \frac{1}{\epsilon N} \sum_{y=1}^{[\epsilon N]} \xi(x+y) \, ,
\end{equation}
which we denote by $W_N(\epsilon,H,\xi)$. By the entropy inequality 
$$
\mathbb{E}_{\mu^{\scriptscriptstyle{N}}}^N \left[ \max_{1\le i\le k} \left\{ \int_0^T ds W_N(\epsilon,H_i(s,\cdot),\xi_s) \right\} \right]
$$
is bounded by
$$
\frac{H(\mu^N|\nu^N_{\alpha , \gamma})}{N} + \frac{1}{N} \log \mathbb{E}_{\nu^{\scriptscriptstyle{N}}_{\alpha ,\gamma}}^N \left[ \exp \left\{ N \max_{1\le i\le k} \int_0^T ds W_N(\epsilon,H_i(s,\cdot),\xi_s)  \right\} \right] .
$$
Now hypothesis (E5) and the elementary inequality $e^{\max a_i} \le \sum e^{a_i}$ imply that this last expression is bounded by
$$
C + \frac{1}{N} \log \mathbb{E}_{\nu^{\scriptscriptstyle{N}}_{\alpha , \gamma}}^N \left[ \sum_{i=1}^k \exp \left\{ N \int_0^T ds W_N(\epsilon,H_i(s,\cdot),\xi_s)  \right\} \right] .
$$ 
Here, since $\max \{ \limsup_{N} N^{-1} \log a_N , \limsup_{N} N^{-1} \log b_N \}$ is greater or equal to $\limsup_{N} N^{-1} \log \{ a_N + b_N \}$, the second term is dominated by
\begin{equation}
\label{eq:W1}
\max_{1\le i\le k} \limsup_{N\ra \oo} \frac{1}{N} \log \mathbb{E}_{\nu^{\scriptscriptstyle{N}}_{\alpha ,\gamma}}^N \left[ \exp \left\{ N \int_0^T ds W_N(\epsilon,H_i(s,\cdot),\xi_s)  \right\} \right] .
\end{equation}
Analogously to the proof of (\ref{eq:derPst}) in Lemma \ref{lemma:replacementbymean}, we have that the previous expression is bounded by
\begin{equation}
\label{eq:W2}
\max_{1\le i\le k} \int_0^T ds \sup_f \left\{ \int W_N(\epsilon,H_i(s,\cdot),\xi)f(\xi) \nu_{\alpha}(d\xi) - N D(f) \right\} + \beta (1-\alpha)^{-1} \, ,
\end{equation}
where the supremum is taken over all densities $f$ with respect to $\nu_{\alpha}$. We just have to estimate the first term in the previous formula, which we are now going to show that it is in fact non-positive. Since, from (\ref{eq:expansao}),
$$
\int \{\xi(x) - \xi(x+[\epsilon N])\} f(\xi) \nu_{\alpha}(d\xi)
= \!\!\!\! \sum_{y=x}^{x+[\epsilon N] -1} \!\!\!\!\! \int \xi(y+1) [1-\xi(y)] \{ f(\xi^{y,y+1}) -f(\xi)\} \nu_{\alpha} (d\xi) , 
$$
we have, for each $B>0$, that
\begin{eqnarray}
\lefteqn{ H(s,x/N) \int  \{\xi(x) - \xi(x+[\epsilon N])\} f(\xi) \nu_{\alpha}(d\xi) \le } \nn \\
& & \le \frac{B}{2} \sum_{y=x}^{x+[\epsilon N] -1} \int \xi(y+1) [1-\xi(y)] \left\{ \sqrt{f(\xi^{y,y+1})} - \sqrt{f(\xi)} \right\}^2 \nu_{\alpha} (d\xi) + \nn \\
& & + \frac{H(s,x/N)^2}{2B} \sum_{y=x}^{x+[\epsilon N] -1} \! \int \xi(y+1) [1-\xi(y)] \left\{ \sqrt{f(\xi^{y,y+1})} + \sqrt{f(\xi)} \right\}^2 \nu_{\alpha} (d\xi) \nn .
\end{eqnarray}
Hence
\begin{eqnarray*}
\lefteqn{ \sum_{x=1}^\oo H(s,x/N) \int \{\xi(x) - \xi(x+[\epsilon N])\} f(\xi) \nu_{\alpha}(d\xi) \le } \nn \\
& & \le \frac{[\epsilon N] B}{2} D(f) + \frac{2}{B} \sum_{x=1}^\oo H(s,x/N)^2 \sum_{y=0}^{[\epsilon N]} \int \xi(x+y) f(\xi) \nu_{\alpha}(d\xi) . \nn 
\end{eqnarray*}
Taking $B=2N$ we obtain from (\ref{eq:W}) that
$$
\int W_N(\epsilon,H(s,\cdot),\xi)f(\xi) \nu_{\alpha}(d\xi) \le N D(f) \, .
$$
Thus the first term in (\ref{eq:W2}) is non-positive and (\ref{eq:W1}) is bounded by $\beta (1-\alpha)^{-1}$. $\square$

\medskip
\no \textbf{Proof of Theorem \ref{theorem:orderofint}:} Let $\mathbb{Q}^*$ be a limit point of the sequence $\mathbb{Q}^N$. By Lemma \ref{lemma:energyestimate} for $\mathbb{Q}^*$ almost every path $\zeta(t,u)$ there exists $B=B(\zeta)>0$ such that
\begin{eqnarray}
\label{eq:energyestimate}
\int_0^T ds \left\{ \int_{\mathbb{R}_+} du \partial_u H(s,u) \zeta(s,u) + H(s,0) \lim_{\epsilon \ra 0} \epsilon^{-1} \int_0^{\epsilon} \zeta(s,u)du \right\} -& &  \nn \\
 -  2 \int_0^T ds \int_{\mathbb{R}_+} du H(s,u)^2 \le B \, , \qquad \qquad & &
\end{eqnarray}
for every $H \in C^{0,1}_K([0,T],\mathbb{R}_+)$. Note that, since $\zeta<1$ we were able to suppress it in the last integrand. Equation (\ref{eq:energyestimate}) implies that 
$$ 
\lambda(H) := \int_0^T ds \left\{ \int_{\mathbb{R}_+} du \partial_u H(s,u) \zeta(s,u) + H(s,0) \lim_{\epsilon \ra 0} \epsilon^{-1} \int_0^{\epsilon} \zeta(s,u)du \right\}
$$
is a bounded linear functional on $C^{0,1}_K([0,T],\mathbb{R}_+)$ for the $L^2$-norm. Since we have that $C^{0,1}_K([0,T],\mathbb{R}_+)$ is a dense subset of $L^2([0,T],\mathbb{R}_+)$, we extend this functional to a bounded linear functional on $L^2([0,T],\mathbb{R}_+)$. By Riesz Representation Theorem, there exists a $L^2$ function, denoted by $\vartheta(s,u)$, such that
$$
\lambda(H)= - \int_0^T ds \int_{\mathbb{R}_+} du H(s,u) \vartheta(s,u) \, ,
$$
for every smooth function $H:[0,T] \times \mathbb{R}_+ \ra \mathbb{R}$ with compact support. Since the integration by parts formula holds for $\vartheta(s,u)$, by definition it is the weak partial derivative of $\zeta(s,u)$ in the second variable. $\square$   

\medskip
\subsection{Removal of the entropy condition}
\label{sec:rementcond}

We have proved Theorem \ref{theorem:hbep} under condition (E5) which is required in our proof of Lemma \ref{lemma:replacementbymean}. If we want to show Theorem \ref{theorem:hbep} without such entropy condition and we follow the steps for the proof described in section \ref{sec:proof}, we see that we only need to show that Lemma \ref{lemma:weaksol} remains true. To obtain Lemma \ref{lemma:weaksol} the idea is to couple the dissipative system not to one but to two exclusion processes such that the exclusion processes are themselves coupled in a appropriate way and for one of them the initial conditions satisfy (E5). This basic coupling for the system, described as in section \ref{sec:system}, is a Feller process on $\Omega_0 = \mathbb{N}^{\mathbb{Z}_{-}} \times \{0,1\}^{\mathbb{Z}_{+}^*} \times \{0,1\}^{\mathbb{Z}_{+}^*}$ whose generator acting on local functions $F: \Omega_0 \ra \mathbb{R}$ is given by
\begin{eqnarray}
\mathcal{L}_0 F(\eta, \xi, \bar{\xi}) & = &
\frac{1}{2} \sum_{|x-y|=1, \, x,y\ge 1 \atop{\xi(x)=\bar{\xi}(x)=1 \atop{ \xi(y)=\bar{\xi}(y)=0} } }  [F(\eta, \xi^{x,y}, \bar{\xi}^{x,y}) - F(\eta, \xi, \bar{\xi})] \nn \\
& & + \frac{1}{2} \sum_{|x-y|=1, \, x,y\ge 1 \atop{\xi(x)=1, \xi(y)=0 \atop {\bar{\xi}(y)=1 \mathrm{ or } \bar{\xi}(x)=0} } }  [F(\eta, \xi^{x,y}, \bar{\xi}) - F(\eta, \xi, \bar{\xi})] \nn \\
& & + \frac{1}{2} \sum_{|x-y|=1, \, x,y\ge 1 \atop{\bar{\xi}(x)=1, \bar{\xi}(y)=0 \atop {\xi(y)=1 \mathrm{ or } \xi(x)=0} } }  [F(\eta, \xi, \bar{\xi}^{x,y}) - F(\eta, \xi, \bar{\xi})] \nn \\
& & + \frac{1}{2} \sum_{|x-y|=1, \, x,y\le 0}  g(\eta(x)) [F(\sigma^{x,y} \eta,\xi,\bar{\xi}) - F(\eta,\xi,\bar{\xi})] \nn \\
& & + g(\eta(0)) [F(\eta - \varrho_0,\tau \xi, \tau \bar{\xi}) - F(\eta,\xi)] \nn \, ,
\end{eqnarray}
where the notation in the above expression is taken from section \ref{sec:system}. For existence results and properties of such a coupling see 
chapter 8 of \cite{liggett}. One important property is that this coupling preserves stochastic order in the sense that: if we consider two inital conditions for the system $\mu_1 = \mu^- \times \mu_1^+$ and $\mu_2 = \mu^- \times \mu_2^+$ such that $\mu_1^+$ is stochastically dominated by $\mu_2^+$ then there exists a coupling measure on $D(\mathbb{R}_+,\Omega_0)$ concentrated on  
$\{ (\eta_t,\xi_t,\bar{\xi}_t) \in \Omega_0 : \xi_t\le \bar{\xi}_t\}$,
for every $0\le t \le T$, and with marginals $\mathbb{P}_{\mu_1}$ with respect to $(\eta,\xi)$ and $\mathbb{P}_{\mu_2}$ with respect to $(\eta,\bar{\xi})$.

The order preserving property may even be described in a less restrictive sense if we improve a bit our coupling considering that the particles are all distinct and that once a $\xi$ particle is attached to a $\bar{\xi}$ particle, which means that from the moment they share the same site, they remain attached from this moment on. This coupling will be called the Stirring coupling for the system. Then the property we have mentioned is the following: consider a set $\Lambda \in \mathbb{Z}_+^*$ and for a measure $\mu$ on $\Omega$ let $\mu^\Lambda$ be the marginal of $\mu$ on $\Lambda$, i.e.,
$$
\mu^{\Lambda} (\varsigma)
= \mu \{ (\eta,\xi) : \xi(x) = \varsigma(x) \textrm{ for all } x \in \Lambda \}, \textrm{ for all } x \in \{0,1\}^{\Lambda} .\nn
$$ 
Fix $\Lambda$ and consider two inital conditions for the system $\mu_1 = \mu^- \times \mu_1^{\Lambda} \times \mu_1^{\mathbb{Z}_+^* - \Lambda}$ and $\mu_2 = \mu^- \times \mu_2^{\Lambda} \times \mu_2^{\mathbb{Z}_+^* - \Lambda}$ such that $\mu_1^{\Lambda}$ is stochastically dominated by $\mu_2^{\Lambda}$. Denote by
$$
\begin{array}{lll}
\mathcal{K}_s^t (\varsigma, \Lambda ,\Gamma) & = &\textrm{number of } \varsigma \textrm{ particles at sites of } \Gamma \textrm{ at time t } \\
& & \textrm{that were at sites of } \Lambda \textrm{ at time s} 
\end{array}
$$
for every trajectory $(\varsigma_t)_{t\ge0}$ on $\{0,1\}^{\mathbb{Z}^*_+}$ and $s<t$.
Thus for the Stirring coupled process with marginals $\mathbb{P}_{\mu_1}$ with respect to $(\eta,\xi)$ and $\mathbb{P}_{\mu_2}$ with respect to $(\eta,\bar{\xi})$, we have that for any other subset $\Gamma$ of $\mathbb{Z}_+^*$, we have that $
\mathcal{K}_0^t (\xi, \Lambda ,\Gamma) \le \mathcal{K}_0^t (\bar{\xi}, \Lambda ,\Gamma)$, 
for almost all trajectories $(\eta,\xi,\bar{\xi})$ on $\Omega_0$ with respect to the coupling measure.

\smallskip

Let $\{\mu^N:N \ge 1\}$ be a sequence in $\mathcal{P}_\pm (\Omega)$ associated to a initial profile $\zeta_0:\mathbb{R}_+ \ra \mathbb{R}$. For each fixed $M$ denote by $\{\mu^{N,M}:N \ge 1\}$ the measure $\mu^{N,-}\times \mu^{N,\Lambda_{N,M}}\times \nu_{1/2}^{\mathbb{Z}_+^*-\Lambda_{N,M}}$, for $\Lambda_{N,M} = \{1,...,NM\}$. The sequence $\{\mu^{N,M}:N \ge 1\}$ is associated to the profile 
$\zeta_0^M(x) := \zeta_0 \mathbf{1}\{[0,M]\}(x) + (1/2) \mathbf{1}\{(M,+\oo)\}(x)$.
Moreover, if $\{\mu^N:N \ge 1\}$ satisfies conditions (E1)-(E4) then $\{\mu^{N,M}:N \ge 1\}$ satisfies conditions (E1)-(E5) and Theorem \ref{theorem:hbep}, in particular Lemma \ref{lemma:weaksol}, holds under this family of initial probability measures. We are going to consider the Stirring coupling between the process speeded up by $N^2$ starting at $\mu^N$ with the same process speeded up by $N^2$ starting at $\mu^{N,M}$. Denote by $(\eta,\xi_t^N, \xi_t^{N,M})$ the coupled process and by $\mathbb{P}^N_{(\mu^{\scriptsize{N}},\mu^{\scriptsize{N,M}})}$ the induced measure on $D([0,T],\Omega_0)$.

\smallskip

It is simple to verify using the next lemma that Lemma \ref{lemma:weaksol} holds under $\{\mu^N:N \ge 1\}$ as the initial family of probability measures

\begin{lemma}
\label{lemma:comparacao}
For every $0 \le c_1 < c_2$ and every continuous function with compact support $H:\mathbb{R}_+ \ra \mathbb{R}$, we have that for all $\delta >0$ 
\begin{equation}
\label{eq:meanMN}
\lim_{M\ra \oo} \lim_{N\ra \oo} \mathbb{P}^N_{(\mu^{\scriptsize{N}},\mu^{\scriptsize{N,M}})} \left[
\sup_{0\le t \le T} \frac{1}{N} \sum_{x=c_1N}^{c_2N} \{\xi_t^{N,M}(x) - \xi_t^N(x)\} > \delta \right] = 0
\end{equation}
and 
\begin{equation}
\label{eq:meanHMN}
\lim_{M\ra \oo} \lim_{N\ra \oo} \mathbb{P}^N_{(\mu^{\scriptsize{N}},\mu^{\scriptsize{N,M}})} \Bigg[
\sup_{0\le t \le T} \frac{1}{N} \sum_{x\ge1} H(x/N)\{\xi_t^{N,M}(x) - \xi_t^N(x)\} > \delta \Bigg] = 0 \, .
\end{equation}
\end{lemma}

\medskip

We conclude this section with the proof of Lemma \ref{lemma:comparacao} which is based on the attractiveness of the system through the order preserving property of the Stirring coupling. 

\medskip  

\no \textbf{Proof:} Fix $0 \le c_1 < c_2$ and let $H:\mathbb{R}_+ \ra \mathbb{R}$ be a continuous function with compact support. Without loss of generality, suppose that the support of $H$ is in $[0,c_2]$. Put $C=c_2$. Note that the supremum in (\ref{eq:meanMN}) is bounded by
$$
\frac{ \mathcal{K}_0^t ( \xi^N, \{MN+1,\dots\}, \{1,\dots,CN\}) + \mathcal{K}_0^t (\xi^{N,M}, \{MN+1,\dots\}, \{1,\dots,CN\})}{N} \, .
$$
and the supremum in (\ref{eq:meanHMN}) is bounded by this same expression times $\|H\|_\oo$. Now, let $\xi$ be an exclusion process with arbitrary initial condition. Then $\mathcal{K}_0^t(\xi,\Gamma, \Lambda)$ is dominated by the same quantity related to the process starting at configurations with a particle on each site in $\Lambda$. Based on this, Let $\nu^{N,M}= \mu^{N,-} \times \nu^{N,M,+}$ be the measure on $\Omega$ such that $\nu^{N,M,+}$ is the Bernoulli product measure on $\mathbb{Z}_+^*$ with marginals given by $\nu^{N,M}\{\xi(x)=1\} = \mathbf{1}\{(NM,+\oo)\} (x)$. Denote by $\bar{\xi}_t^{M,N}$ the process starting with the measure $\nu^{N,M}$. Consider the Stirring coupling of the system with starting measures $\mu^N$ and $\nu^{N,M}$ and then with starting measures $\mu^{N,M}$ and $\nu^{N,M}$. From the property of the coupling, we have that
$$
\max [\mathcal{K}_0^t (\xi^N, \{MN+1,\dots\}, \{1,\dots,CN\}) , \mathcal{K}_0^t (\xi^{N,M}, \{MN+1,\dots\}, \{1,\dots,CN\}) ] 
$$
is bounded by
$\mathcal{K}_0^t (\bar{\xi}^{N,M}, \{MN+1,\dots\}, \{1,\dots,CN\})$, which is equal to the number of particles at sites $\{x\in \mathbb{Z} : x\le CN\}$ at time t, for the process starting at $\nu^{N,M}$. This number is clearly bounded by the number of particles at sites 
in $(-\oo, CN]$ at time t for the simple symmetric excusion process (construct a coupling similar to the previous one on $\mathbb{Z}$), which divided by $N$ converges to the integral on $(-\oo ,C]$ of the solution of
$$
\left\{
\begin{array}{l}
\partial_t \zeta(t,u) = \frac{1}{2} \Delta \zeta(t,u), \quad t \in \mathbb{R}_+, \ u \in \mathbb{R}, \\
\zeta(0,u)= \mathbf{1}\{(M,+\oo)\}(u), \quad u \in \mathbb{R} \, .
\end{array}
\right. 
$$
By Diffusion Theory, the unique solution of the previous equation has a stochastic representation given by $\mathbb{E}_u \left[ h \! ( B_t ) \right], \ u\in\mathbb{R} , \ t\in\mathbb{R}_+,$
where $(B_t)$ is a standard Brownian Motion and $h=\mathbf{1}\{(M,+\oo)\}$. It is then a straightforward computation, using the Gaussian kernel, to show that the integral on $(-\oo ,C]$ of this last expression converges to 0 exponentially fast as $M\ra \oo$. Thus (\ref{eq:meanMN}) holds. $\square$

\bigskip
\section{From the exclusion process to the Potts model}
\setcounter{equation}{0}
\label{sec:exctopotts}

We prove in this section Theorem \ref{theorem:potts}. As a first identification we associate to each configuration $f$ in $\mathcal{I}$ such that $f(0)=0$ a configuration in $\mathbb{N}^{\mathbb{Z}}$ representing the increments of the former: $\eta(x)=f(x+1)-f(x)$ for every $x \in \mathbb{Z}$. This allows us to associate the Potts model dynamics to a zero range dynamics as described in section \ref{sec:potts}.

For technical reasons we consider the zero range as two coupled processes: the dissipative and the absorbing systems. A configuration $\eta$ in $\mathbb{N}^{\mathbb{Z}}$ is associated to a configuration $(\eta,\xi)$ in $\Omega=\mathbb{N}^{\mathbb{Z}_-} \times \{0,1\}^{\mathbb{Z}_+^*}$ in such a way that, for $x \ge 1$, $\eta(x)$ represents the number of consecutive holes that precede the $x^{th}$ particle in configuration $\xi$. Since for the exclusion process the total number of sites in a given finite box equals the total number of holes plus the total number of particles, we obtain the following relation:
$$
\sum_{x=1}^{n} \xi(x) + \sum_{y=1}^{\sum_{x=1}^{n} \xi(x)} \eta(y) \le n
\le \sum_{x=1}^{n} \xi(x) + \sum_{y=1}^{1 +\sum_{x=1}^{n} \xi(x)} \eta(y) \, ,
$$
for all $n\ge 1$
This is the same as
\begin{equation}
\label{eq:eptozr}
\frac{1}{N} \sum_{y=1}^{N \left\{ \frac{1}{N} \sum_{x=1}^{AN} \xi(x) \right\}} \eta(y) \le A - \frac{1}{N} \sum_{x=1}^{AN} \xi(x) \le \frac{1}{N} \sum_{y=1}^{N \left\{ \frac{1}{N} + \frac{1}{N} \sum_{x=1}^{AN} \xi(x) \right\} } \eta(y)
\end{equation}
for every $A>0$. To a probability measure $\mu$ on $\mathbb{N}^{\mathbb{Z}}$ let $\mathcal{V}\mu$ be the probability measure on $\Omega$ that corresponds to the push-forward of $\mu$ through the map described above. Now, Fix a sequence $\{ \mu^N:N\ge 1\}$ of probability measures on $\mathbb{N}^{\mathbb{Z}}$ such that $\{ \mathcal{V}\mu^N:n\ge 1\}$ is a sequence of probability measures that satisfies the conditions of Theorem \ref{theorem:hbep}. Let $\zeta: \mathbb{R}_+ \times \mathbb{R}_+ \ra \mathbb{R}$ be the unique solution of (\ref{eq:pdeep}) with initial condition equals to the initial density profile for $\{ \mathcal{V}\mu^N:n\ge 1\}$. Then from Theorem \ref{theorem:hbep} and inequality (\ref{eq:eptozr}), for every $B>0$, the mean $N^{-1} \sum_{y=0}^{BN} \eta_{tN^2}(y)$ converges in probability to $M^{-1}(t,B) - B$, where
$$
M(t,A) = \int_0^A \zeta(t,u) du ,
$$
which is an increasing absolutely continuous function on the second variable. Writing
$$
\rho(t,u) = \partial_u (M^{-1}(t,u) - u), \ u > 0,
$$
we obtain from the definition of $M$ that
$$
\rho(t,u)=\frac{1}{\zeta(t,M^{-1}(t,u))}-1,  \ u > 0.
$$
This function $\rho: \mathbb{R}_+ \times \mathbb{R}_+ \ra \mathbb{R}$ is the unique solution of (\ref{eq:pdezr}) with initial condition
$$
\rho_0(u)=\frac{1}{\zeta_0(M^{-1}(0,u))}-1,  \ u > 0.
$$
Thus we have proved:

\begin{theorem}
\label{theorem:hbzr}
Fix a sequence of $\{\mu^N:N\ge 1\}$ on $\mathbb{N}^{\mathbb{Z}}$ such that $\{ \mathcal{V}\mu^N:n\ge 1\}$ is a sequence of probability measures that satisfies the conditions of Theorem \ref{theorem:hbep}. Then, for every continuous function $G:\mathbb{R} \ra \mathbb{R}$ 
$$
\lim_{N\ra {\oo}} \mathbb{P}_{\mu^{\scriptscriptstyle{N}}}^N \left[ \left| \frac{1}{N} \sum_{x \in \mathbb{Z}} G(x/N)\eta_t(x) - \int duG(u)\rho(t,u) \right| \ge \delta \right] = 0
$$
for every $0\le t\le T$ and $\delta>0$, where $\rho$ is the unique solution of (\ref{eq:pdezr}).
\end{theorem}

Now fix a smooth function $G:\mathbb{R}\ra\mathbb{R}$ with compact support. Note that
$$
\frac{1}{N} \sum_{x \in \mathbb{Z}} G(x/N) N^{-1} [f_{tN^2}(x) - f_{tN^2}(0)] = \frac{1}{N} \sum_{x \in \mathbb{Z}} \left\{ \frac{1}{N} \sum_{y\ge x+1} G(y/N) \right\} \eta_t(x).
$$
From theorem \ref{theorem:hbzr} the term at the right of this equation converges in probability to
$$
\int_{\mathbb{R}_+} du \, G(u) \lambda(t,u),
\quad \textrm{where} \quad
\lambda(t,u) = \int_0^u \rho(t,v)dv.
$$
Therefore, $\lambda$ is the unique solution of \ref{eq:pdepm}. This concludes the proof of Theorem \ref{theorem:potts}.

\bigskip

\no \emph{Acknowledgements.} The author wish to thank Claudio Landim, his advisor during the prepararion of this paper.

\bigskip


\end{document}